\documentclass[12pt,a4paper]{article}
\usepackage{latexsym,epsfig,fullpage,amsfonts,amssymb,amsmath,amscd,epic,graphics,
rotating,theorem,yhmath,pxfonts,comment,tikz}
\usepackage[all]{xy}
\theoremstyle{plain}
\newtheorem{lemma}{Lemma}[section]
\newtheorem{proposition}[lemma]{Proposition}
\newtheorem{remark}[lemma]{Remark}

\newtheorem{theorem}[lemma]{Theorem}

\newtheorem{corollary}[lemma]{Corollary}
\theorembodyfont{\rmfamily}  \font\ncsc=cmcsc10
 \font\ntt=cmtt12

\tikzstyle{regular_node}=[circle, draw, fill=black!50,
                        inner sep=0pt, minimum width=4pt]
\tikzstyle{small_node}=[circle, draw, fill=black!50,
                        inner sep=0pt, minimum width=2pt]

\begin{document}
\newcommand{\pperp}{\hbox{$\perp\hskip-6pt\perp$}}
\newcommand{\ssim}{\hbox{$\hskip-2pt\sim$}}
\newcommand{\aleq}{{\ \stackrel{3}{\le}\ }}
\newcommand{\ageq}{{\ \stackrel{3}{\ge}\ }}
\newcommand{\aeq}{{\ \stackrel{3}{=}\ }}
\newcommand{\bleq}{{\ \stackrel{n}{\le}\ }}
\newcommand{\bgeq}{{\ \stackrel{n}{\ge}\ }}
\newcommand{\beq}{{\ \stackrel{n}{=}\ }}
\newcommand{\cleq}{{\ \stackrel{2}{\le}\ }}
\newcommand{\cgeq}{{\ \stackrel{2}{\ge}\ }}
\newcommand{\ceq}{{\ \stackrel{2}{=}\ }}
\newcommand{\A}{{\mathbb A}}
\newcommand{\K}{{\mathbb K}}
\newcommand{\Z}{{\mathbb Z}}\newcommand{\F}{{\mathbf F}}
\newcommand{\R}{{\mathbb R}}
\newcommand{\C}{{\mathbb C}}
\newcommand{\Q}{{\mathbb Q}}
\newcommand{\PP}{{\mathbb P}}
\newcommand{\mnote}{\marginpar}\newcommand{\red}{{\operatorname{red}}}
\newcommand{\Leaf}{{\operatorname{Leaf}}}\newcommand{\Sym}{{\operatorname{Sym}}}
\newcommand{\oeps}{{\overline\eps}}\newcommand{\Div}{{\operatorname{Div}}}
\newcommand{\oDel}{{\widetilde\Del}}\newcommand{\lab}{{\operatorname{lab}}}
\newcommand{\real}{{\operatorname{Re}}}\newcommand{\Aut}{{\operatorname{Aut}}}
\newcommand{\conv}{{\operatorname{conv}}}\newcommand{\BG}{{\operatorname{BG}}}
\newcommand{\Span}{{\operatorname{Span}}}\newcommand{\GS}{{\operatorname{GS}}}
\newcommand{\Ker}{{\operatorname{Ker}}}\newcommand{\Ss}{{\operatorname{SS}}}
\newcommand{\Ann}{{\operatorname{Ann}}}\newcommand{\mt}{{\operatorname{wt}}}
\newcommand{\Fix}{{\operatorname{Fix}}}\newcommand{\Ima}{{\operatorname{Im}}}
\newcommand{\sign}{{\operatorname{sign}}}\newcommand{\Def}{{\operatorname{Def}}}
\newcommand{\Tors}{{\operatorname{Tors}}}
\newcommand{\Card}{{\operatorname{Card}}}
\newcommand{\alg}{{\operatorname{alg}}}\newcommand{\ord}{{\operatorname{ord}}}
\newcommand{\oi}{{\overline i}}
\newcommand{\oj}{{\overline j}}
\newcommand{\bb}{\boldsymbol{b}}
\newcommand{\os}{{\overline s}}
\newcommand{\ba}{\boldsymbol{a}}
\newcommand{\be}{\boldsymbol{e}}
\newcommand{\ow}{{\overline w}}
\newcommand{\bc}{\boldsymbol{c}}
\newcommand{\oz}{{\overline z}}\newcommand{\on}{{\overline n}}
\newcommand{\eps}{{\varepsilon}}
\newcommand{\proofend}{\hfill$\Box$\bigskip}
\newcommand{\Int}{{\operatorname{Int}}}\newcommand{\grad}{{\operatorname{grad}}}
\newcommand{\pr}{{\operatorname{pr}}}
\newcommand{\Hom}{{\operatorname{Hom}}}
\newcommand{\Ev}{{\operatorname{Ev}}}
\newcommand{\im}{{\operatorname{Im}}}\newcommand{\br}{{\operatorname{br}}}
\newcommand{\sk}{{\operatorname{sk}}}\newcommand{\Fl}{{\operatorname{Fl}}}
\newcommand{\const}{{\operatorname{const}}}
\newcommand{\Sing}{{\operatorname{Sing}}\hskip0.06cm}
\newcommand{\conj}{{\operatorname{Conj}}}
\newcommand{\Cl}{{\operatorname{Cl}}}
\newcommand{\Crit}{{\operatorname{Crit}}}
\newcommand{\Ch}{{\operatorname{Ch}}}
\newcommand{\discr}{{\operatorname{discr}}}
\newcommand{\Tor}{{\operatorname{Tor}}}
\newcommand{\Conj}{{\operatorname{Conj}}}
\newcommand{\vol}{{\operatorname{vol}}}
\newcommand{\defect}{{\operatorname{def}}}
\newcommand{\codim}{{\operatorname{codim}}}
\newcommand{\tmu}{{\C\mu}}
\newcommand{\bv}{\boldsymbol{v}}
\newcommand{\ox}{{\overline{x}}}
\newcommand{\bw}{{\boldsymbol w}}\newcommand{\bn}{{\boldsymbol n}}
\newcommand{\bx}{{\boldsymbol x}}
\newcommand{\bd}{{\boldsymbol d}}
\newcommand{\bz}{{\boldsymbol z}}\newcommand{\bp}{{\boldsymbol p}}
\newcommand{\tet}{{\theta}}
\newcommand{\Del}{{\Delta}}
\newcommand{\bet}{{\beta}}
\newcommand{\kap}{{\kappa}}
\newcommand{\del}{{\delta}}
\newcommand{\sig}{{\sigma}}
\newcommand{\alp}{{\alpha}}
\newcommand{\Sig}{{\Sigma}}
\newcommand{\Gam}{{\Gamma}}
\newcommand{\gam}{{\gamma}}\newcommand{\idim}{{\operatorname{idim}}}
\newcommand{\Lam}{{\Lambda}}
\newcommand{\lam}{{\lambda}}
\newcommand{\SC}{{SC}}
\newcommand{\MC}{{MC}}
\newcommand{\nek}{{,...,}}
\newcommand{\cim}{{c_{\mbox{\rm im}}}}
\newcommand{\clM}{\tilde{M}}
\newcommand{\clV}{\bar{V}}

\title{Refined descendant invariants of toric surfaces}
\author{Lev Blechman
\and Eugenii Shustin}
\date{}
\maketitle
\begin{abstract}
We construct refined tropical enumerative genus zero invariants of toric surfaces that specialize to the tropical descendant genus zero invariants introduced by Markwig and Rau when the quantum parameter
tends to $1$. In the case of trivalent tropical curves our invariants turn to be the
G\"ottsche-Schroeter refined broccoli invariants. We show that this is the only possible refinement of the Markwig-Rau descendant invariants that generalizes the G\"ottsche-Schroeter refined broccoli invariants. We discuss also the computational aspect
(a lattice path algorithm) and exhibit some examples.
\end{abstract}

{\it MSC-2010}: 14N10, 14T05

\medskip


\section*{Introduction}

Starting from Mikhalkin's foundational work \cite{Mi} tropical geometry has served as an ultimate tool to solve important enumerative problems. Later it has become clear that tropical geometry provides new insights to various problems of ``classical" geometry. The present work has been inspired by two such phenomena. One is that the genus zero descendant invariants of the plane and other toric surfaces, defined as integrals over the moduli spaces of stable maps of rational curves, can be computed via enumeration of certain plane tropical curves \cite{MR}. Another exciting phenomenon is the existence of refined tropical enumerative invariants, i.e., tropical enumerative invariants depending on a parameter \cite{BG,GS,IM}. We will comment on this in more detail.

F. Block and L. G\"ottsche \cite{BG} defined
a refined multiplicity of a plane trivalent tropical curve
being a symmetric Laurent polynomial
in one variable $y^{1/2}$. They showed that, for $y=1$, its value is the Mikhalkin
multiplicity as introduced in \cite[Definition
4.15]{Mi}) and, for $y=-1$, its value is the Welschinger multiplicity as introduced in
\cite[Definition 7.19]{Mi}. Notice that, under appropriate conditions, enumeration
of plane trivalent tropical curves with Mikhalkin and
Welschinger multiplicities gives Gromov-Witten and Welschinger invariants of toric del Pezzo surfaces, respectively
(see \cite[Theorems 1 and 6]{Mi}). In general, the enumerative meaning of the refined count
remains in question, though in certain cases it is related to
the quantization of real plane curves in the sense of
\cite{Mi1}.

I. Itenberg and G. Mikhalkin \cite{IM} showed that the total refined multiplicity of the
plane tropical curves
having given degree and genus and passing through an appropriate number of generic
points in the plane
does not depend on the choice of the point constraints. This is a generalization of the invariance of
the count of tropical curves of given degree and genus with Mikhalkin and Welschinger
multiplicities established earlier in
\cite[Theorem 4.8]{GM} and \cite[Theorem 1]{IKS}.

L. G\"ottsche and F. Schroeter \cite{GS} advanced further and introduced another refined tropical invariant of genus zero,
which, for $y=-1$ specializes to the so-called broccoli invariant (see \cite{GMS}) and, for $y=1$ specializes to
the genus zero tropical descendant invariant $\big\langle\tau_0(2)^{n_0}
\tau_1(2)^{n_1}\big\rangle_\Del$ of a toric surface associated to a convex lattice polygon $\Del$ as introduced in \cite{MR}.

{\bf The main result} of this paper is a construction of a refinement of an arbitrary genus zero tropical descendant invariant
$\big\langle\prod_{k\ge0}\tau_k(2)^{n_k}\big\rangle_\Del$; we present such
a refined invariant in
Theorem \ref{tpsi1}, Section \ref{spsi1}. The construction is inductive with the G\"ottsche-Schroeter invariants as the base, while in the induction step, we
replace the vertices of valency $>3$ with trees having vertices of smaller valency.
In Section \ref{secu}, we show that the combinatorial type of trees
used in the induction step does not matter,
and one
always obtains the same refined invariant, which generalizes the Block-G\"ottsche and G\"ottsche-Schroeter invariants.

Similarly to the Block-G\"ottsche and G\"ottsche-Schroeter invariants
our invariant counts certain rational plane tropical curves through a generic configuration of points, and each
curve is counted with a refined multiplicity equal to the product of multiplicities of vertices
(normalized by the automorphisms).
This property allows one to use combinatorial computational tools like the
``lattice path algorithm". In Appendix to this paper,
we indicate how to modify the Markwig-Rau lattice path algorithm \cite[Section 9]{MR} in order to
compute our refined invariant and provide an example of computation.

Another refined tropical count related to descendant Gromov-Witten invariants has been suggested in
\cite{M}. The tropical objects counted there are plane tropical disks (i.e., halves of tropical curves) which have
only trivalent vertices equipped with the refined Block-G\"ottsche multiplicities. In turn, we consider the
entire tropical curves admitting multi-valent marked vertices, and the main novelty of our work is
the definition of refined multiplicities of marked vertices of any valency.

We remark that our invariant is not in general a Laurent polynomial in $y$ (or even in $y^{1/2}$)
if we allow marked vertices
of valency $>3$, but it may
have poles at $y=0$ and $y=-1$. We estimate the order of the pole at $y=-1$ in
Proposition \ref{lpsi4} (Section \ref{spsi1}) and show that our upper bound is sharp
(Section \ref{sec-ex}).

An interesting question is how to extend the definition of tropical
descendant invariants and their refinements to curves of positive genus.
The particular case of genus one curves with at most trivalent marked vertices
has been settled in
\cite{SS}. Another perspective question is whether the refined tropical descendant invariant can be interpreted via Mikhalkin's quantization of algebraic curves \cite{Mi1}.

\section{Plane marked rational tropical curves}
We shortly recall some basic definitions concerning plane rational tropical curves, adapted to our setting (for details, see \cite{GM,GM1,MR,Mi2,Mi,NS} and especially \cite{GKM}).

\subsection{Plane rational tropical curves}\label{sec-tc}
Under an {\it $n$-marked plane rational tropical curve} we understand
a triple $(\Gam,\bp,h)$, where
\begin{itemize}\item $\Gam$ is a finite
connected metric tree without vertices of valency $\le2$, whose set $\Gam^0$ of vertices in nonempty, and the
set of edges $\Gam^1$ contains a
subset $\Gam^1_\infty\ne\emptyset$ consisting of edges isometric to $[0,\infty)$ (called {\it ends}), while
$\Gam^1\setminus\Gam^1_\infty$ consists of edges isometric to compact segments in $\R$ (called {\it finite edges});
\item $h:\Gam\to\R^2$ is a proper continuous map such that
$h$ is nonconstant, affine-integral on each edge of $\Gam$ in the length coordinate
(i.e., given by $(x,y)=(a,b)t+(x_0,y_0)$ with $(a,b)\in\Z^2\setminus\{0\}$ and $t$
being the length coordinate on the considered edge)
and, at each vertex $V$ of $\Gam$, the {\it balancing condition} holds
\begin{equation}\sum_{E\in\Gam^1,\ V\in E}
\ba_V(E)=0\ ,\label{e1}\end{equation} where $\ba_V(E)$
is the image under the differential $D(h\big|_E)$
of the unit tangent vector to $E$ emanating from
its endpoint $V$ (called the {\it directing vector} of $E$ centered at $V$);
\item $\bp=(p_1,...,p_n)$ is a sequence of $n$ points of $\Gam$.
\end{itemize}

We call a vertex $V\in\Gam^0$ {\it flat} if the vectors in the
left-hand side of (\ref{e1}) span a one-dimensional subspace of $\R^2$

Notice that each vector $\ba_V(E)$ can
be written as $\ba_V(E)=m\bv$, where $m$ is a positive integer
(called the {\it weight} of the edge $E$ of $(\Gam,\bp,h)$)
and $\bv\in\Z^2\setminus\{0\}$ is primitive.
The {\it degree} of the plane rational tropical curve
$(\Gam,\bp,h)$ is the multi-set of vectors
$$\Del(\Gam,\bp,h)=\left\{\ba_V(E)),\ E\in\Gam^1_\infty\right\}\ .$$
The balancing condition yields that $\Del(\Gam,\bp,h)$ is a {\it balanced}
multi-set, i.e.
$$\sum_{\bb\in\Del(\Gam,\bp,h)}\bb=0\ .$$
We call $\Del(\Gam,\bp,h)$ {\it primitive} if it contains only primitive
vectors, and we call
$\Del(\Gam,\bp,h)$ {\it nondegenerate}, if $\Span_\R\Del(\Gam,\bp,h)=\R^2$.

With a plane rational tropical curve $(\Gam,h)$ we associate
an {\it unparameterized} (or {\it embedded}) {\it plane tropical curve} in the sense of
\cite[Section 2.1]{Mi2} or \cite[Section 2]{Mi}, which we denote
$T=h_*(\Gam)$. This is a closed finite rational\footnote{Here
``rational" means ``with rational slopes".} one-dimensional
polyhedral subcomplex of $\R^2$ supported at $h(\Gam)$.
Its edges are equipped with
positive integral weights as follows: for an edge $e$ of $T$ pick any of its interior points
$x\in e$, and define the weight $w(e)$ of $e$ as the sum of the weights of those edges of $\Gam$, whose $h$-images cover $x$. It is easy to see that the balancing condition (\ref{e1}) yields
the balancing condition for $T$ at each vertex $v\in T^0$ in the form
$$\sum_{v\in e} w(e)\cdot \ba_v(e)=0\ ,$$ where $e$ ranges over all the edges of $T$ incident to $v$,
and $\ba_v(e)$ denotes the primitive integral vector directing $e$ and emanating from $v$.
The degree $\Del(T)$ of $T$ is the multi-set of vectors $w(e)\cdot \ba_v(e)$, where $e$ runs over all unbounded edges of $T$.
If we take the vectors $\bb\in\Del(T)$ in the natural cyclic order, rotate them by $\frac{\pi}{2}$ clockwise,
and attach the initial point of each vector to the endpoint of the preceding one, we obtain a convex lattice
polygon $P(T)$ (denoted also $P(\Del)$), called the {\it Newton polygon} of $T$ and of $(\Gam,\bp,h)$. There is
a duality
between the edges and vertices of $T$ on one side and the
edges and
polygons of a certain (dual) subdivision of $P(T)$, see \cite[Proposition 2.1]{Mi2}. We denote the dual object by
${\mathcal D}(*)$. The weight $w(e)$ of an edge $e$ of $T$ equals the lattice length
of the dual edge ${\mathcal D}(e)$ in the above subdivision.

By a {\it combinatorial type} of a plane rational marked tropical curve $(\Gam,\bp,h)$ we mean
the combinatorial type of the pair $(\Gam,\bp)$ enhanced with the collection of
directing vectors $\ba_V(E)$ for all edges $E\in\Gam^1$ and vertices $V\in E$.

\subsection{Moduli spaces of plane rational marked tropical curves}
The moduli spaces of plane rational marked tropical curves
are our main objects of study. Here we recall some information on moduli spaces,
following \cite{GKM,GM,GM1,MR,Mi3,Mi} and adapting notations to our setting.

Under an isomorphism $(\Gam,\bp,h)\to(\Gam',
\bp',h')$ we understand an isometry
$\varphi:\Gam\to\Gam'$, identifying the ordered sequences $\varphi:\bp\to\bp'$, and
satisfying $h=h'\circ\varphi$.
Isomorphism classes
of plane rational $n$-marked tropical curves of degree $\Del$ are parameterized
by the moduli space ${\mathcal M}_{0,n}(\R^2,\Del)$.

We will also use labeled tropical curves. In this case, we fix a linear order on $\Del$ and denote the obtained sequence by $\Del^\lab$. A {\it labeled $n$-marked plane rational tropical curve of degree}
$\Del^\lab$ is a triple $(\Gam^\lab,\bp,h)$, where $(\Gam,\bp,h)$ is an $n$-marked plane
rational tropical curve of degree $\Del$, and $\Gam^\lab$ is the graph $\Gam$ with a linear order on the set of its ends $\Gam^1_\infty$ such that $h$ equips the $i$-th end of $\Gam^\lab$ with the
$i$-th directing vector in $\Del^\lab$ for all $i$
(cf. \cite[Definition 4.1]{GKM}).
Notice that a labeled metric tree $\Gam^\lab$ does not admit nontrivial automorphisms. Thus, the
moduli space ${\mathcal M}^\lab_{0,n}(\R^2,\Del^\lab)$ parameterizes just $n$-marked plane rational tropical curves of degree $\Del^\lab$. According to \cite[Section 4, Lemma 4.6 and Proposition 4.7]{GKM} (see also \cite[Lemma 2.1]{GM1}), the geometry of ${\mathcal M}_{0,n}(\R^2,\Del)$
and ${\mathcal M}^\lab_{0,n}(\R^2,\Del^\lab)$ can be described as follows.

\begin{lemma}\label{sl1}
(1) The space ${\mathcal M}^\lab_{0,n}(\R^2,\Del^\lab)$ can be identified with a polyhedral fan
of pure dimension $|\Del|-1+n$ in some Euclidean space $\R^N$. Open cells of this fan are in bijection with the combinatorial types of the labeled curves
$(\Gam^\lab,\bp,h)\in{\mathcal M}^\lab_{0,n}(\R^2,\Del^\lab)$, while independent parameters are:
the coordinates of $h(V)\in\R^2$ for a chosen vertex $V\in\Gam^\lab$, the lengths of the finite edges whose interior is disjoint from $\bp$, and the distances from each marked point $p_i\in\bp\setminus\Gam^0$ to an endpoint of the edge $E\supset\{p_i\}$. The faces of a cell correspond to the case of vanishing of some parameters. The top-dimensional cells correspond to $\Gam$ trivalent and $\bp\cap\Gam^0=\emptyset$.

(2) The group $G$ of permutations $\varphi:\Del^\lab\to\Del^\lab$ such that $\varphi(\bb)=\bb$ for each $\bb\in\Del^\lab$ acts on ${\mathcal M}^\lab_{0,n}(\R^2,\Del^\lab)$, and one has the finite surjective quotient map
$$\pi_{0,n}:{\mathcal M}^\lab_{0,n}(\R^2,\Del^\lab)\to{\mathcal M}_{0,n}(\R^2,\Del)\ .$$
Furthermore, for any element $[(\Gam,\bp,h)]\in{\mathcal M}_{0,n}(\R^2,\Del)$, we have
\begin{equation}\big|\pi_{0,n}^{-1}([(\Gam,\bp,h)])\big|=\frac{|G|}{|\Aut(\Gam,\bp,h)|}
\ .\label{elab}\end{equation}
\end{lemma}

Furthermore, for any sequence $\on=(n_k)_{k\ge0}\in\Z_+^\infty$, where $\Z_+=\{m\in\Z\ :\ m\ge0\}$ and $n=\sum_{k\ge0}n_k$, introduce
\begin{eqnarray}{\mathcal M}_{0,\on}(\R^2,\Del)&=&\Big\{(\Gam,\bp,h)\in{\mathcal M}_{0,n}(\R^2,\Del)\ :\ \nonumber\\
& & p_i\in\bp\ \text{are interior points of edges for}\ 1\le i\le n_0,\nonumber\\
& & p_i\in\bp\ \text{are}\ (k+2)-\text{valent vertices for}\ \sum_{j<k}n_j<i\le\sum_{j\le k}n_j,\quad k\ge 1\Big\}\ .
\nonumber\end{eqnarray}
Denote by $\widehat{\mathcal M}_{0,\on}(\R^2,\Del)$ the closure
of ${\mathcal M}_{0,\on}(\R^2,\Del)$ in ${\mathcal M}_{0,n}(\R^2,\Del)$. Respectively, we have the labeled analogue ${\mathcal M}^\lab_{0,\on}(\R^2,\Del^\lab)\subset\widehat{\mathcal M}^\lab_{0,\on}(\R^2,\Del^\lab)
\subset{\mathcal M}^\lab_{0,n}(\R^2,\Del^\lab)$.
The following statement is straightforward from Lemma \ref{sl1}.

\begin{lemma}\label{l5lab}
The space $\widehat{\mathcal M}^{\;\lab}_{0,\on}(\R^2,\Del^\lab)$
is either empty, or is
a finite polyhedral fan of pure dimension $|\Del|-1+\sum_{k\ge0}(1-k)n_k$.
Its open top-dimensional cells parameterize
tropical curves $(\Gam^\lab,\bp,h)$ with $\Gam\setminus\bp$ trivalent and
$\bp\cap\Gam^0$ consisting of exactly $\sum_{i\ge1}n_i$ points among which $n_i$ points are $(i+2)$-valent vertices of
$\Gam$ for all $i\ge1$. The space $\widehat{\mathcal M}_{0,\on}(\R^2,\Del)$
is the quotient of $\widehat{\mathcal M}^{\;\lab}_{0,\on}(\R^2,\Del^\lab)$ by the action of the group $G$.
\end{lemma}

Introduce the evaluation maps
$$\Ev:\widehat{\mathcal M}^{\;\lab}_{0,\on}(\R^2,\Del^\lab)\to\widehat{\mathcal M}_{0,\on}(\R^2,\Del)\to
\R^{2n},\quad
\Ev(\Gam,\bp,h)=h(\bp)\ .$$
Suppose that
a sequence
$\on=(n_k)_{k\ge0}\in\Z_+^\infty$ satisfies $\dim\widehat{\mathcal M}^{\;\lab}_{0,\on}(\R^2,\Del^\lab)=2n$, that is,
\begin{equation}\sum_{k\ge0}(k+1)n_k=|\Del|-1\label{epsi1},\quad\sum_{k\ge0}n_k=n\ .\end{equation}
Denote by ${\mathcal M}^{e,\lab}_{0,\on}(\R^2,\Del^\lab)$, resp.,
${\mathcal M}^{e}_{0,\on}(\R^2,\Del)$,
the union of those open cells of dimension
$2n=|\Del|-1+\sum_{k\ge0}(1-k)n_k$ in ${\mathcal M}^{\lab}_{0,\on}(\R^2,\Del^\lab)$,
resp., ${\mathcal M}_{0,\on}(\R^2,\Del)$, whose $\Ev$-images have dimension $2n$, and by
$\widehat{\mathcal M}^{\;e,\lab}_{0,\on}(\R^2,\Del^\lab)$, resp.,
$\widehat{\mathcal M}^{\;e}_{0,\on}(\R^2,\Del)$, the closure in ${\mathcal M}^\lab_{0,n}(\R^2,\Del^\lab)$, resp., $\widehat{\mathcal M}_{0,n}(\R^2,\Del)$.
By $\Ev^e$
we denote the
restriction of $\Ev$ to $\widehat{\mathcal M}^{\;e,\lab}_{0,\on}(\R^2,\Del^\lab)$ and
$\widehat{\mathcal M}^{\;e}_{0,\on}(\R^2,\Del)$.

Let us recall the regularity notion due to Mikhalkin \cite[Sections 2.6 and 4.5]{Mi}.
We say thata curve $(\Gam,\bp,h)$ is {\it regular}, if each component $K$
of the set $\Gam\setminus\bp$
is {\it regular}, that is, $K$ is trivalent, has no flat vertices and contains exactly one unbounded edge.
The following statement is a reformulation  of \cite[Lemmas 4.20 and 4.22]{Mi}.

\begin{lemma}\label{lreg}
Let a component $K$ of $\Gam\setminus\bp$ for some plane rational tropical curve
$(\Gam,\bp,h)$ be regular. Then the following holds.

(1) The edges of the closure $\overline K\subset\Gam$ admit a unique orientation (called the regular
orientation) such that the marked points are sources, the unbounded edges are oriented towards infinity, and
each vertex of $K$ is incident to exactly two incoming edges.

(2) Denote by $(\Gam_K,\bp_K,h_K)$ the plane rational tropical curve, where $\bp_K=\bp\cap\overline K$, and $\Gam_K$ is obtained
from $\overline K$ by extending each edge incident to a marked point to an unbounded edge,
while $h_K$ linearly extends $h$ beyond the marked points. Put $\Del_K=\deg(\Gam_K,\bp_K,h_K)$ and $n'=|\bp_K|$. Then the map
$$\Ev:{\mathcal M}_{0,n'}(\R^2,\Del_K)\to\R^{2n'}$$ defines a linear isomorphism of the germ\footnote{Here and further on, under the {\it germ} we understand
a sufficiently small Euclidean neighborhood of the central element.} of
${\mathcal M}_{0,n'}(\R^2,\Del_K)$ at $[(\Gam_K,\bp_K,h_K)]$ onto a germ of $\R^{2n'}$ at $h_K(\bp_K)$.
\end{lemma}

Lemma \ref{lreg} has an interesting consequence for us:

\begin{lemma}\label{lne}
For an arbitrary nondegenerate balanced multiset $\Del\subset\Z^2\setminus\{0\}$ and a
sequence $\on=(n_i)_{i\ge0}\in\Z^\infty_+$ satisfying (\ref{epsi1}), the moduli space
${\mathcal M}^e_{0,\on}(\R^2,\Del)$ contains an element represented by
a regular curve $(\Gam,\bp,h)$. Moreover, the induced embedded plane tropical curve
$T=h_*(\Gam)$ is dual to a subdivision of the Newton polygon $P(\Del)$ into $n_0+2\sum_{i\ge1}n_i-1$
nondegenerate convex polygons, obtained by drawing $n_0+2\sum_{i\ge1}n_i-2$ chords that
join certain pairs of integral points in $\partial P(\Del)$ and do not intersect each other in their interior points.
\end{lemma}

It follows that the moduli spaces ${\mathcal M}^e_{0,\on}(\R^2,\Del)\subset{\mathcal M}_{0,\on}(\R^2,\Del)$ are nonempty, which
strengthens the second assertion of Lemma \ref{l5lab}.

{\bf Proof of Lemma \ref{lne}.} Suppose that there exists a regular element $[(\Gam,\bp,h)]\in{\mathcal M}_{0,\on}(\R^2,\Del)$.
Then, by Lemma \ref{lreg}, it belongs to ${\mathcal M}^e_{0,\on}(\R^2,\Del)$.
Furthermore, denote by $m_1$ the number
of the unmarked trivalent vertices of such a curve $(\Gam,\bp,h)$.
Comparing the Euler characteristic $\chi(\Gam)=1$, the number $|\Gam^1|$, and formula
(\ref{epsi1}), we obtain
\begin{equation}m_1=\sum_{i\ge0}n_i-1\ .\label{etrival}\end{equation}

We construct the curve $(\Gam,\bp,h)$ asserted in the lemma using the ideas of
\cite[Proposition 2.10]{IM}, where the case of $\on=(n_0,0,0,...)$ was considered.

The boundary of the Newton polygon $P(\Del)$ can be represented as the union of
cyclically ordered integral segments $[v_k,v_{k+1}]$, $k=0,...,|\Del|-1$,
$v_{|\Del|}=v_0$, obtained by
rotating the vectors $\ba_k\in\Del$, $k=0,...,|\Del|-1$, by $\frac{\pi}{2}$ clockwise (see Section \ref{sec-tc}).
The set ${\mathcal V}=\{v_0,...,v_{|\Del|-1}\}$ includes all the vertices of $P(\Del)$.

We prove the lemma by induction on $\sum_{i\ge1}n_i$. For the base of induction, suppose that $n_i=0$ for all $i\ge1$.
Then there exists a subdivision of $P(\Del)$ into $m_1$ (nondegenerate) triangles obtained by drawing
$m_1-1$ chords joining some pairs of points in ${\mathcal V}$ so that no two chords intersect
in their interior points. This is easily derived by induction on $|\Del|$: if $|\Del|\ge4$, then there exist
two non-collinear segments $[v_{i-1},v_i]$, $[v_i,v_{i+1}]$ such that $[v_{i-1},v_{i+1}]\not\subset\partial\Del$.
The obtained triangulation is convex, i.e., lifts to a graph of a convex
piecewise-linear function $\nu:P(\Del)\to\R$, and hence, the subdivision is dual to a plane tropical curve
$T$, which is an embedded trivalent tree. Picking a marked point on all but one ends of $T$, we obtain a regular tropical curve as desired.

\begin{figure}
\setlength{\unitlength}{1cm}
\begin{picture}(14,5)(0,0)
\thinlines
\put(1,3){\line(1,0){4}}

\dashline{0.2}(7,2)(8,2)\dashline{0.2}(9.5,1.5)(10,1)

\thicklines
\put(1,2){\line(0,1){1}}\put(1,2){\line(2,-1){2}}
\put(3,1){\line(2,1){2}}\put(5,2){\line(0,1){1}}
\put(1,3){\line(4,-1){4}}\put(1,3){\line(1,1){1}}
\put(2,4){\line(2,1){2}}\put(5,3){\line(-1,2){1}}

\put(8.1,2){\line(1,0){0.9}}\put(9,2){\line(1,-1){0.5}}
\put(9,2){\line(1,3){0.6}}\put(9.6,3.8){\line(1,0){1}}
\put(9.6,3.8){\line(0,1){0.7}}\put(9.6,4.5){\line(2,1){0.6}}
\put(9.6,4.5){\line(-1,1){0.7}}
\put(9.6,4.5){\line(-1,2){0.4}}

\put(2.8,0){(a)}\put(8.8,0){(b)}
\put(9.2,2.7){$\ast$}\put(9.5,4.4){$\bullet$}
\put(0.5,2.9){$v_i$}\put(5.2,2.9){$v_{i+k+1}$}\put(5.2,1.9){$v_{i+k+2}$}
\put(2.9,1.8){$P'$}\put(2.9,3.5){$P''$}
\put(11,1.7){$(\Gam',\bp',h')$}\put(11,3.7){$(\Gam'',\bp'',h'')$}

\end{picture}
\caption{Proof of Lemma \ref{lne}}\label{fig-div}
\end{figure}
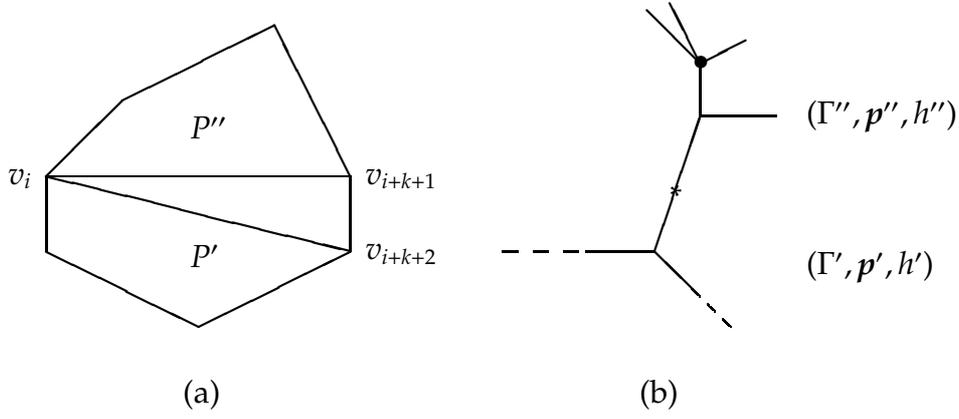

For the induction step suppose that $n_k>0$ for some $k\ge1$. If $|\Del|=k+2$, then $\Gam$ has one $(k+2)$-valent marked vertex
and $|\Del|=k+2$ ends incident to it. The tropical curve $T=h_*(\Gam)$ is dual to the entire polygon
$P(\Del)$. If $|\Del|=k+3$, then formulas (\ref{epsi1}) and (\ref{etrival}) imply that $n_0=n_k=1$, $n_i=0$, $i\ne0,k$, and
$m_1=1$. As noticed in the preceding paragraph,
there exist two non-collinear segments $[v_{i-1},v_i]$, $[v_i,v_{i+1}]$ such that $[v_{i-1},v_{i+1}]\not\subset\partial\Del$.
The chord $[v_{i-1},v_{i+1}]$ defines the required subdivision of $P(\Del)$, while
the corresponding curve $(\Gam,\bp,h)$ has a marked $(k+2)$-valent vertex joined by a bounded edge orthogonal to
the above chord with an unmarked trivalent vertex, and an extra marked point is chosen in one of the ends incident to
the unmarked vertex. In the remaining case $|\Del|\ge k+4$, we have $m_1\ge1$ by formulas (\ref{epsi1}) and (\ref{etrival}).
We claim that there exists $0\le i\le|\Del|-1$ such that the chord $[v_i,v_{i+k+1}]$ is not contained in $\partial P(\Del)$.
Indeed, otherwise, either all points $v_i,...,v_{i+k+1}$ lie on one line, or all points
$v_{i+k+1},...,v_i$ lie on one line (here we follow the cyclic order in ${\mathcal V}$), and this holds for all $i=0,...,|\Del|-1$,
which is only possible when all points $v_0,...,v_{|\Del|-1}$ lie of one line, contrary to the assumption that
$\Del$ is nondegenerate. For a similar reason, at least one of the chords
$[v_i,v_{i+k+2}]$, $[v_{i-1},v_{i+k+1}]$ is not contained in $\partial P(\Del)$. If, for instance,
$[v_i,v_{i+k+2}]\not\subset\partial P(\Del)$, then we obtain the polygons (see Figure \ref{fig-div}(a))
\begin{itemize}\item $P'=\conv(v_{i+k+2},v_{i+k+3},...,v_i)=P(\Del')$, where the
nondegenerate, balanced multiset $\Del'$
consists of the vectors $\ba_j\in\Del$, $j\not\in\{i,...,i+k+1\}$, and the vector
$\ba'=\ba_i+...+\ba_{i+k+1}$;
\item $P''=\conv(v_i,v_{i+1},...,v_{i+k+2})=P(\Del'')$, where the nondegenerate, balanced
multiset $\Del''$ consists of the vectors $\ba_i,...,\ba_{i+k+1}\in\Del$ and the
vector $-\ba'$.
\end{itemize} Let $\on'=(n'_i)_{i\ge0}$, $n'_k=n_k-1$, $n'_i=n_i$, $i\ne k$. Since $\Del'$ and $\on'$
satisfy (\ref{epsi1}), and $\sum n'_i<\sum n_i$, by induction assumption we get
a curve $(\Gam',\bp',h')\in{\mathcal M}^e_{0,\on'}(\R^2,\Del')$ matching the requirements of Lemma.
We also get a plane tropical curve $(\Gam'',\bp'',h'')$ of degree $\Del''$,
having a marked $(k+2)$-valent vertex and an unmarked trivalent vertex and associated with the subdivision of
$P''$ by the chord $[v_i,v_{i+k+1}]$.
Finally, we obtain the desired element $[(\Gam,\bp,h)]\in{\mathcal M}^e_{0,\on}(\R^2,\Del)$
by gluing the curves $(\Gam',\bp',h')$ and $(\Gam'',\bp'',h'')$ along their ends directed by the vectors $\ba'$ and $-\ba'$, respectively (see Figure \ref{fig-div}(b)).
\proofend

The next statement provides a geometric background for the proof of the invariance of the refined enumeration
of tropical curves introduced in the present paper.
It is very similar to considerations in \cite[Sections 3 and 4]{MR}.

\begin{proposition}\label{l1} Let $\Del\subset\Z^2\setminus\{0\}$ be a non-degenerate balanced multi-set, and let
a sequence $\on=(n_k)_{k\ge0}\in\Z_+^\infty$ satisfy (\ref{epsi1}). Then
the target space of the map \mbox{$\Ev^e:\widehat{\mathcal M}^{\:e,\lab}_{0,\on}(\R^2,\Del^\lab)\to\R^{2n}$}
splits into the disjoint union
$$\R^{2n}=X^{2n}\cup X^{2n-1}\cup X^{2n-2}\ ,$$ where
\begin{enumerate}\item[(1)] $X^{2n}$ is the union of finitely many open polyhedra of dimension $2n$, and,
for each element $\bx\in X^{2n}$, its preimage $(\Ev^e)^{-1}(\bx)$ is
finite; furthermore, each curve $(\Gam^\lab,\bp,h)\in (\Ev^e)^{-1}(\bx)$ belongs to
${\mathcal M}^\lab_{0,\on}(\R^2,\Del^\lab)$ and is regular;
\item[(2)] $X^{2n-1}$ is the union of finitely many (relatively) open polyhedra of dimension $2n-1$, and, for each point $\bx\in X^{2n-1}$, its preimage $(\Ev^e)^{-1}(\bx)$ is
finite;
furthermore, each curve $(\Gam^\lab,\bp,h)\in(\Ev^e)^{-1}(\bx)$ is as follows:
\begin{enumerate}\item[(2i)] either $(\Gam,\bp,h)$ is as in item (1);
\item[(2ii)] or $(\Gam^\lab,\bp,h)\in{\mathcal M}^\lab_{0,\on}(\R^2,\Del^\lab)$, and all but one components of
the set $\Gam\setminus\bp$ are regular, while the remaining component has one unbounded edge, one four-valent vertex, and
the rest of vertices are trivalent;
\item[(2iii)] or $(\Gam^\lab,\bp,h)\in{\mathcal M}^\lab_{0,\on'}(\R^2,\Del^\lab)$, where $\on'=(n'_i)_{i\ge0}\in\Z_+^\infty$,
and there exists $k\ge0$ such that $n_k>0$, $n'_k=n_k-1$, $n'_{k+1}=n_{k+1}+1$, and $n'_i=n_i$ for all $i\ne k,k+1$;
furthermore, all but one components of the set $\Gam\setminus\bp$ are regular, while the remaining component is
bounded, trivalent, and incident to one of the $(k+3)$-valent marked vertices;
\end{enumerate}
\item[(3)] $X^{2n-2}$ is a closed finite polyhedral complex of dimension $\le2n-2$.
\end{enumerate}
\end{proposition}

{\bf Proof.}
Define
$$X^{2n}=\R^{2n}\setminus\Ev^e\left(\widehat{\mathcal M}^{\;e,\lab}_{0,\on}(\R^2,\Del^\lab)
\setminus{\mathcal M}^{e,\lab}_{0,\on}(\R^2,\Del^\lab)\right)\ .$$ Now, let ${\mathcal M}^{2n-1}$ be the union
of those $(2n-1)$-dimensional cells of $\widehat{\mathcal M}^{\;e,\lab}_{0,\on}(\R^2,\Del^\lab)
\setminus{\mathcal M}^{e,\lab}_{0,\on}(\R^2,\Del^\lab)$, which are injectively projected into $\R^{2n}$ by $\Ev^e$. Then define
$$X^{2n-1}=\Ev^e({\mathcal M}^{2n-1})\setminus\Ev^e\left(\widehat{\mathcal M}^{\;e,\lab}_{0,\on}(\R^2,\Del^\lab)
\setminus({\mathcal M}^{e,\lab}_{0,\on}(\R^2,\Del^\lab)\cup{\mathcal M}^{2n-1})\right)\ ,$$
$$X^{2n-2}=\Ev^e\left(\widehat{\mathcal M}^{\;e,\lab}_{0,\on}(\R^2,\Del^\lab)
\setminus({\mathcal M}^{e,\lab}_{0,\on}(\R^2,\Del^\lab)\cup{\mathcal M}^{2n-1})\right)\ .$$
So, claim (3) follows by construction.

The finiteness of $(\Ev^e)^{-1}(\bx)$ for each $\bx\in X^{2n}\cup X^{2n-1}$ follows from Lemma
\ref{l5lab} and the fact that the cells of ${\mathcal M}^{e,\lab}_{0,\on}(\R^2,\Del^\lab)\cup{\mathcal M}^{2n-1}$ are injectively mapped
into $\R^{2n}$.

Let us show the regularity of any curve $(\Gam^\lab,\bp,h)\in(\Ev^e)^{-1}(\bx)$, $\bx \in X^{2n}$. Note that
no component $\gamma$ of $\Gam\setminus\bp$ is bounded. Indeed, otherwise, it would yield a constraint to the position of
the images of the marked points on $\partial\gamma$, that is, the image of the corresponding $2n$-cell of
${\mathcal M}^e_{0,\on}(\R^2,\Del)$ would be of dimension $<2n$ (cf. \cite[Lemma 4.20]{Mi}).
Since the number of connected component of $\Gam\setminus\bp$ is $1+\sum_{k\ge0}(k+1)n_k=|\Del|$
(see (\ref{epsi1})), we derive from the above observation that each component of $\Gam\setminus\bp$ contains
exactly one unbounded edge. Notice that there are no flat unmarked vertices: indeed, otherwise,
contrary to the finiteness of $(\Ev^e)^{-1}(\bx)$ we would obtain a one-parameter family inside
$(\Ev^e)^{-1}(\bx)$ when varying the position of such a flat vertex along the line containing the images of its incident edges.
Finally, all components of $\Gam\setminus\bp$ must be trivalent due to the condition of maximal-dimensional image
(cf. \cite[Proposition 2.23]{Mi}).

For claim (2), note that by construction, elements of any $(2n-1)$-cell of $\widehat{\mathcal M}^{\;e,\lab}_{0,\on}(\R^2,\Del^\lab)$ are degenerations of elements of some $2n$-cell. Hence, they appear via either moving one of the marked points outside
$\Gam^0$ to a vertex of $\Gam$, or contracting exactly one bounded edge. In the former case, we fit the situation (2iii) with $k=0$.
In the latter case, either a marked vertex collates with an unmarked, trivalent one, which fits the situation (2iii) with some $k>0$,
or two unmarked trivalent vertices collate, that is, conditions of (2ii) are satisfied. Note also that, in item (2ii), the regular orientation
on the complement to the marked points for elements of the $2n$-cell that degenerate to the considered curve $(\Gam^\lab,\bp,h)$, induces an orientation
on the edges of the special component of $\Gam\setminus\bp$ such that three edges incident to the unmarked four-valent vertex
are incoming and one is outgoing.
\proofend

\section{Refined count of plane rational marked tropical curves}\label{sec1}

Throughout this section, we fix a standard basis in $\R^2$, and
for any $\ba=(a_1,a_2)$, $\bb=(b_1,b_2)\in\R^2$,
set $\ba\wedge\bb=\det\left(\begin{matrix}a_1&a_2\\ b_1&b_2\end{matrix}
\right)$. We also set
\begin{equation}[\alpha]_y^-=\frac{y^{\alpha/2}-y^{-\alpha/2}}{y^{1/2}-y^{-1/2}},\quad
[\alpha]_y^+=\frac{y^{\alpha/2}+y^{-\alpha/2}}{y^{1/2}+y^{-1/2}},\quad\text{for all}\ \alpha\in\R\ ,
\label{ed2}\end{equation} $y$ being a formal parameter.

\subsection{Refined multiplicity of a plane rational marked tropical curve}\label{spsi2}
Let us be given a nondegenerate, balanced multi-set $\Del\subset\Z^2\setminus\{0\}$ and its linearly ordered form
$\Del^\lab$, a positive integer $n$ and
a sequence $\on=(n_k)_{k\ge0}\in\Z_+^\infty$ satisfying (\ref{epsi1}). Let $(\Gam,\bp,h)\in
{\mathcal M}^e_{0,\on}(\R^2,\Del)$ be regular, and let
$(\Gam^{\;\lab},\bp,h)\in
{\mathcal M}^{e,\lab}_{0,\on}(\R^2,\Del)$ be one of the labeled preimages of
$(\Gam,\bp,h)$.
We start with defining a refined multiplicity $RM_y(\Gam,\bp,h,V)$ (depending on
a formal parameter $y$) for each vertex $V\in\Gam^0$.

\smallskip

\noindent {\it (1) Refined multiplicity of a trivalent vertex.}
Suppose that $V\in\Gam^0$ is trivalent and choose two distinct edges $E_1,E_2\in\Gam^1$ incident
to $V$. Define the {\it Mikhalkin's multiplicity} of the vertex $V$ by
(cf. \cite[Definition 2.16]{Mi})
\begin{equation}\mu(\Gam,h,V)=|\ba_V(E_1)\wedge\ba_V(E_2)|
\ .\label{epsi3}\end{equation} Due
to the balancing condition
(\ref{e1}), this number does not depend on the choice of the pair of edges incident to $V$
and, in fact,
is equal to the lattice area of the triangle ${\mathcal D}(h(V))$, dual to the vertex $h(V)$
of the embedded plane
tropical curve $T=h_*(\Gam)$. Following \cite{BG,GS}, we define
\begin{equation}
RM_y(\Gam,\bp,h,V)=\begin{cases}[\mu(\Gam,h,V)]_y^+\ ,\quad &\text{if}\ V\in\bp,\\
[\mu(\Gam,h,V)]_y^-\ ,\quad &\text{if}\ V\not\in\bp.\end{cases}
\label{epsi4}\end{equation}

\smallskip

{\it (2) Refined multiplicity of a marked vertex of valency $\ge4$.} For any balanced sequence
$$A=(\ba_i)_{i=1,...,m},\quad m\ge2,\quad\ba_i\in\R^2,\ i=1,...,m,\quad\sum_{i=1}^m\ba_i=0\ ,$$
we will recursively define an expression $\theta_y(A)$ containing a formal parameter $y$.

If $m=2$, we set $\theta_y(A)=1$. If $m=3$, we set
$\theta_y(A)=[|\ba_1\wedge\ba_2|]^+_y$. Note that, due to the balancing condition, this definition does not depend
on the choice of order in the sequence $A$. Furthermore, it holds that
$\theta_y(A)=RM_y(\Gam,\bp,h,V)$ when $V\in\Gam^0\cap\bp$ is trivalent, and $\ba_i=\ba_V(E_i)$,
$i=1,2,3$, with $E_1,E_2,E_3$ being the edges of $\Gam$ incident to $V$. 

If $m\ge4$, then, for each pair $1\le i<j\le m$, we form the two balanced sequences
\begin{itemize}\item $A'_{ij}$ consisting of the vectors $\ba_k$, $1\le k\le m$, $k\ne i,j$, and one more vector
$\ba_{ij}:=\ba_i+\ba_j$,
\item $A''_{ij}=(\ba_i,\ba_j,-\ba_{ij})$.
\end{itemize} Then we set
\begin{equation}\theta_y(A)=\sum_{1\le i<j\le m}\theta_y(A'_{ij})\cdot\theta_y(A''_{ij})\ .\label{emu}\end{equation}
It is easy to see that $\theta_y(A)$ does not depend on the choice of the order in $A$.
At last, observe that $\theta_y(A)$ can be written as the sum over all
plane rational trivalent tropical curves of degree $A$ counted with multiplicity proportional to the
product of the factors $[\mu(V)]^+_y$ over all trivalent vertices $V$ of a given curve.

Now, given a vertex $V\in\Gam^0\cap\bp$ of valency $m$ and somehow ordered edges $E_1,...,E_m$ of $\Gam$ incident to $V$, we define
\begin{equation}RM_y(\Gam,\bp,h,V)=\theta_y(\Del_V),\quad\Del_V=(\ba_V(E_i))_{i=1,...,m}\ .\label{epsi2}\end{equation}

Finally, put
\begin{equation}
RM_y(\Gam^{\;\lab},\bp,h)=
\prod_{V\in\Gam^0}RM_y(\Gam,\bp,h,V),
\quad RM_y(\Gam,\bp,h)=\frac{RM_y(\Gam^{\;\lab},\bp,h)}
{|\Aut(\Gam,\bp,h)|}\ .\label{epsi5}\end{equation}

\subsection{Invariance of the refined count}\label{spsi1}

\begin{theorem}\label{tpsi1}
Let $\Del\subset\Z^2\setminus\{0\}$ be a nondegenerate balanced multi-set, a sequence
$\on=(n_k)_{k\ge0}\in\Z_+^\infty$ satisfy (\ref{epsi1}), and the set
$X^{2n}\subset\R^{2n}$ be as in Proposition
\ref{l1}. Then the expression
\begin{equation}RD_y(\Del,\on,\bx)=\sum_{\renewcommand{\arraystretch}{0.6}
\begin{array}{c}
\scriptstyle{(\Gam,\bp,h)\in{\mathcal M}^e_{0,\on}(\R^2,\Del)}\\
\scriptstyle{h(\bp)=\bx}
\end{array}}RM_y(\Gam,\bp,h)\label{ne-inv}\end{equation}
does not depend on the choice of $\bx\in X^{2n}$.
\end{theorem}

The proof is presented in Sections \ref{sec-plan} and \ref{sec-collision}.

\begin{remark}\label{rpsi1}
Note that the invariant $RD_y$ is an extension of other known refined invariants:
\begin{itemize}
\item if $\on=(n_0,0,...)$, i.e., $n_k=0$ for all $k\ge1$, then $RD_y(\Del,\on)$ coincides with
the genus zero Block-G\"ottsche refined invariant \cite{BG};
\item if $\on=(n_0,n_1,0,...)$, i.e., $n_k=0$ for all $k\ge2$, and all the vectors in $\Del$ are primitive, then
$RD_y(\Del,\on)$ coincides with the G\"ottsche-Schroeter refined invariant \cite{GS}.
\end{itemize}
\end{remark}

The invariant $RD_y(\Del,\on)$ provides a refinement of the tropical
descendant invariants as introduced in \cite{MR}. Namely, introduce the {\it normalized} refined descendant invariant
$$NRD_y(\Del,\on)=\prod_{k\ge1}\left[\frac{3\cdot 2^{k+1}}{(k+2)!(k+1)!}
\right]^{n_k}\cdot RD_y(\Del,\on)\ .$$

\begin{lemma}\label{lpsi2}
Given a nondegenerate balanced multi-set $\Del\subset\Z^2\setminus\{0\}$ and a sequence $\on=(n_k)_{k\ge0}\in\Z_+^\infty$
satisfying (\ref{epsi1}), we have
\begin{equation}\lim_{y\to 1}NRD_y(\Del,\on)=\left\langle\prod_{k\ge0}
\tau_k(2)^{n_k}\right\rangle_\Del\ .\label{epsi11}\end{equation}
\end{lemma}

{\bf Proof.} Along \cite[Theorem 8.4]{MR}, the tropical descendant invariant $\left\langle\prod_{k\ge0}
\tau_k(2)^{n_k}\right\rangle_\Del$ can be computed by counting rational marked tropical curves
in $(\Ev^e)^{-1}(\bx)\subset{\mathcal M}^e_{0,\on}(\R^2,\Del)$ with the multiplicities
$$\omega(\Gam,\bp,h)=\frac{1}{|\Aut(\Gam,\bp,h)|}
\prod_{V\in\Gam^0\setminus\bp}\mu(\Gam,h,V)\ .$$
On the other hand, for a trivalent vertex $V$ of $\Gam$ with $\mu=\mu(\Gam,h,V)$, we
have
$$\lim_{y\to1}RM_y(\Gam,\bp,h,V)=\lim_{y\to1}\frac{y^{\mu/2}-y^{-\mu/2}}{y^{1/2}-y^{-1/2}}=
\mu=\mu(\Gam,h,V)\ ,$$ if $V$ is unmarked, and
$$\lim_{y\to1}RM_y(\Gam,\bp,h,V)=\lim_{y\to1}\frac{y^{\mu/2}+y^{-\mu/2}}{y^{1/2}+y^{-1/2}}=1\ ,$$
if $V$ is marked.
If $V$ is a four-valent marked vertex, then
formula (\ref{emu}) yields six summands, each one equal to $1$, and hence by
(\ref{epsi2}) $RM_1(\Gam,\bp,h,V)=6=4!3!/(3\cdot2^3)$. Then we inductively apply formula
(\ref{emu}) and obtain for any marked vertex of valency $k+2$
$$\lim_{y\to1}RM_y(\Gam^\lab,\bp,h,V)=\frac{(k+2)!(k+1)!}{3\cdot 2^{k+1}},\quad k\ge 2\ .$$ Thus, (\ref{epsi11}) follows.
\proofend

The invariant $RD_y(\Del,\on)$ is often a rational function of $y$:

\begin{proposition}\label{lpsi4}
Given a nondegenerate balanced multi-set
$\Del\subset\Z^2\setminus2\Z^2$ and a sequence $\on=(n_k)_{k\ge0}\in\Z^\infty_+$ satisfying
(\ref{epsi1}), we have
\begin{equation}RD_y(\Del,\on)=\frac{F(y+y^{-1})}{(y+2+y^{-1})^m}\ ,\label{epsi20}\end{equation}
where $m\ge0$ and $F$ is a (nonzero) polynomial of degree
$$\deg F=p_a(P(\Del))+\frac{|\partial P(\Del)\cap\Z^2|-|\Del|}{2}+m\ ,$$ where $p_a(P(\Del))$ is the number of interior integral points
of the Newton polygon $P(\Del)$. Furthermore,
\begin{equation}m\le\sum_{k\ge1}k(n_{2k}+n_{2k+1})\ .
\label{ed1}\end{equation}
\end{proposition}

{\bf Proof.}
To show that $RD_y(\Del,\on)$ is a (rational) function of $y$, we move $y$ around the circle $S_c=\{|y|=c\}$, $0<c\ll1$, and check that
the value of $RD_y(\Del,\on)$ does not change sign. Indeed, notice that the recursion (\ref{emu}) yields that
$RD_y(\Del,\on)$ can be expressed as the sum of multiplicities of finitely many trivalent rational
tropical curves $(\Gam_{(3)},h_{(3)})$
of degree $\Del$, and each multiplicity is of the form $s\prod_{V\in\Gam_{(3)}^0}[\mu(\Gam_{(3)},h_{(3)},V)]_y^{\pm}$ with some $s\in\Q$.
The factor $[\mu(\Gam_{(3)},h_{(3)},V)]_y^{\pm}$ changes its sign as $y$ travels along $S_c$ iff $\mu(\Gam_{(3)},h_{(3)},V)$ is even.
Thus, the claim comes from the statement of \cite[Proposition 2.3(4)]{IM}: if $\Del$ does not contain even vectors, then the number of vertices of
$\Gam_{(3)}$ with even $\mu(\Gam_{(3)},h_{(3)},V)$ is even.

To compute the denominator of the function $RD_y(\Del,\on)$, we consider
a summand $\Sigma=s\prod_{V\in\Gam_{(3)}^0}[\mu(\Gam_{(3)},h_{(3)},V)]_y^{\pm}$ that appears in the
expression of the preceding paragraph. Note that $y^{\mu/2}-y^{-\mu/2}$ is always divisible by $y^{1/2}-y^{-1/2}$ and, in
addition, is divisible by $y^{1/2}+y^{-1/2}$ as $\mu$ is even, and $y^{\mu/2}+y^{-\mu/2}$ is divisible by
$y^{1/2}+y^{-1/2}$ if $\mu$ is odd. Hence, either the denominator in $\Sigma$ cancels out, or is equal to
$(y^{1/2}+y^{-1/2})^{m'-m''}$, where $m'$, resp. $m''$, is the number of marked, resp. unmarked, vertices $V\in\Gam^0_{(3)}$
with even value $\mu(\Gam_{(3)},h_{(3)},V)$. As observed in the preceding paragraph, the number $m'+m''$ is even, and hence the denominator of
$\Sigma$ takes form $(y^{1/2}+y^{-1/2})^{2m}=(y+2+y^{-1})^m$ with an integer $m$.

Formula for $\deg F$ is similar to that in \cite[Proposition 2.10]{IM}. For a summand
$\Sigma=s\prod_{V\in\Gam_{(3)}^0}[\mu(\Gam_{(3)},h_{(3)},V)]_y^{\pm}$ as above, the difference between
the top exponents of $y$ in the nominator and denominator equals
$$\frac{1}{2}\sum_{V\in\Gam_{(3)}^0}(\mu(\Gam_{(3)},h_{(3)},V)-1)=\frac{1}{2}
\sum_{V\in\Gam_{(3)}^0}\mu(\Gam_{(3)},h_{(3)},V)
-\frac{|\Del|}{2}+1\ .$$ Taking into account the geometric meaning of $\mu(\Gam_{(3)},h_{(3)},V)$, we obtain that
the latter expression takes its maximal value for the summands associated with the curve $(\Gam,\bp,h)$
from Lemma \ref{lne}, for which the trivalent trees are obtained by a further subdivision of $P(\Del)$ by chords into
$|\Del|-3$ triangles. Hence, the considered maximal value is
$$\text{Area}(P(\Del))-\frac{|\Del|}{2}+1\quad\overset{\text{Pick's formula}}{=}
\quad p_a(P(\Del))+\frac{|\partial P(\Del)\cap\Z^2|-|\Del|}{2}\ .$$

To establish the bound (\ref{ed1}), we again consider the formula $RD_y(\Del,\on)=\sum_{(\Gam_{(3)},h_{(3)})}
s\prod_{V\in\Gam_{(3)}^0}[\mu(\Gam_{(3)},h_{(3)},V)]_y^{\pm}$, representing the invariant
via the sum over trivalent curves $(\Gam_{(3)},h_{(3)})$.
The number of factors of type $[\mu(\Gam_{(3)},h_{(3)},V)]_y^{+}$ then equals $\sum_{k\ge1}kn_k$, which
yields $(y^{1/2}+y^{-1/2})^{\sum_{k\ge1}kn_k}$ in the denominator. We shall show that
$(y^{1/2}+y^{-1/2})^{\sum_{k\ge1}n_{2k+1}}$ divides the nominator, and thus derive (\ref{ed1}).

Consider a marked vertex of odd valency and a trivalent tree that appears in the computation of the refined multiplicity of that vertex via recursion
(\ref{emu}). The following holds: if Mikhalkin multiplicities of all the trivalent vertices in the tree
are even, then the original marked vertex is incident to an edge of an even weight.
We leave this claim as an easy exercise for the reader.
Following the ideas from \cite[Proof of Proposition 4.2]{SS} and \cite[Proof of Lemma 3.12]{GS},
we introduce a subgraph $\Gam^{even}\subset
\Gam$ including all the edges of even weight and their endpoints.
Now, if $V$ is a marked vertex of odd valency, which does not belong to $\Gam^{even}$, then in each trivalent tree associated with
$V$, there is a vertex of odd Mikhalkin multiplicity. On the other hand, the closure of
each component $K$ of $\Gam^{even}\setminus\bp$ contains at least as many unmarked vertices
of $\Gam$ as the marked ones. Indeed, such a closure $\overline K$ is a tree with univalent and unmarked trivalent vertices. By the Euler characteristic formula, the number of trivalent vertices equals the number of univalent vertices minus $2$. By the regularity condition, there is an unmarked univalent vertex, and hence the number of unmarked
vertices is at least the number of marked ones. Notice that the unmarked vertices of $\overline K$ are, in fact, unmarked trivalent vertices of $\Gamma$ of even Mikhalkin multiplicities (the edges of $\overline K$ all are finite, since the ends of $\Gamma$ have odd weights by the hypotheses of the lemma). It follows that each of the latter vertices provides a factor $y^{1/2}+y^{-1/2}$ in the nominator, which finally yields
at least $\sum_{k\ge1}n_{2k+1}$ factors $y^{1/2}+y^{-1/2}$ in the nominator of $RD_y(\Del,\on)$, and we are done.
\proofend

In general, the denominator in formula (\ref{epsi20}) is unavoidable, and the bound (\ref{ed1}) is sometimes sharp as one can see in the examples of Corollary \ref{ac1} in Section \ref{sec-ex}.

\subsection{Proof of the invariance: Preliminaries}\label{sec-plan}
It will be convenient to consider labeled tropical curves.
In view of formulas (\ref{elab}) and (\ref{epsi5}), the invariance of $RD_y(\Del,\on,\bx)$ is equivalent to the
invariance of $RD^{\lab}_y(\Del,\on,\bx)$.

So, we choose two generic configurations $\bx(0),\bx(1)\in X^{2n}
\subset\R^{2n}$. By the dimension reason, there exists a continuous path $\bx(t)\in\R^{2n}$, $0\le t\le 1$, that joins the chosen configurations, avoids
$X^{2n-2}$, but may finitely many times hit $(2n-1)$-dimensional cells of
$X^{2n-1}$, which may cause changes in the structure of $(\Ev^e)^{-1}(\bx(t))$. We shall consider all possible wall-crossing phenomena
and verify the constancy of
$RD^{\lab}_y(\Del,\on,\bx(t))$ (as a function of $t$) in the in these events.
To relax notations we simply denote labeled tropical curves by $\Gam$ or
$\Gam$ and write $RD^{\lab}_y(\bx(t))$ for $RD^{\lab}_y(\Del,\on,\bx(t))$.

Let $\bx(t^*)$ be generic in an $(2n-1)$-dimensional cell of $X^{2n-1}$. Denote by $H_0$
the germ of this cell at $\bx(t^*)$ and by $H_+,H_-\subset\R^{2n}$ the germs of the halfspaces
with common boundary $H_0$. Let $C^*=(\Gam,\bp,h)\in(\Ev^e)^{-1}(\bx(t^*))$
be as described in Proposition \ref{l1}(2), and let $F_0\subset\widehat{\mathcal M}^{\;e,\lab}_{0,\on}(\R^2,\Del)$
be the germ at $C^*$ of the $(2n-1)$-cell projected by $\Ev^e$
onto $H_0$. We shall analyze the $2n$-cells of $\widehat{\mathcal M}^{\;e,\lab}_{0,\on}(\R^2,\Del)$
attached to $F_0$, their projections onto $H_+,H_-$, and prove the
constancy of $RD^{\lab}_y(\bx(t))$, $t\in(t^*-\eps,t^*+\eps)$, $0<\eps\ll1$.

In the case of Proposition \ref{l1}(2i), the curve $C^*$ deforms keeping its combinatorics as
$t\in(t^*-\eps,t^*+\eps)$, and hence this deformation does not affect the contribution
to $RD^{\lab}_y(\bx(t))$.

In the case of Proposition \ref{l1}(2ii), in the deformation of the curve $C^*$
when $t\in(t^*-\eps,t^*+\eps)$ the multiplicative contributions of the marked and the unmarked trivalent vertices
to the refined multiplicity of the current curve do not change, and one has only to study the deformation
of a neighborhood of the unmarked four-valent vertex. The proof of the constancy of $RD^{\lab}_y(\bx(t))$ in this situation can be found in
\cite[Theorem 1 and Pages 5314--5316]{IM} or in \cite[Section 4.2]{SS}.

\subsection{Proof of the invariance: Collision of a marked and an unmarked vertices}\label{sec-collision}
The remaining task is to show the invariance when crossing the wall described in Proposition
\ref{l1}(2iii).

{\it (1) Preparation.} Let $C^*=(\Gam,\bp,h)\in(\Ev^e)^{-1}(\bx(t^*))$,
where $\bx(t^*)\in X^{2n-1}$,
satisfy the conditions of Proposition \ref{l1}(2iii). That is $\Gam$ has a $(k+3)$-valent vertex
$V\in\Gam^0\cap\bp$, and $C^*$ is the limit of one or several families (which we call $C^*$-families)
$(\Gam(t),\bp(t),h(t))\in(\Ev^e)^{-1}(\bx(t))$,
where either $t\in(t^*,t^*+\eps)$, or $t\in(t^*-\eps,t^*)$, and such that the edge joining some $(k+2)$-valent vertex
$V(t)\in\Gam(t)^0\cap\bp(t)$ and a trivalent vertex $W(t)\in\Gam(t)^0\setminus\bp(t)$ collapses as $t\to t^*$, while
$\lim_{t\to t^*}h(t)\big(V(t)\big)=\lim_{t\to t^*}h(t)\big(W(t)\big)=h(V)$.

Observe that the orientation of the edges of $\Gam$ incident to $V$, which is induced by the regular orientation
of a $C^*$-family, is as follows: one edge (which we denote $E_0$) is oriented towards $V$, and
the remaining edges (which we denote $E_1,...,E_{k+2}$) are oriented outwards; moreover, the edge $E_0$ is the same for all
$C^*$-families, since it is distinguished by the property to be a part of
the unique bounded component of the complement $\Gam\setminus\bp$. Note also, that in each $C^*$-family
$(\Gam(t),\bp(t),h(t))$, the corresponding edge $E_0(t)$ is incident to the unmarked trivalent vertex $W(t)$.

Without loss of generality, we can suppose that, for each $C^*$-family, the image $h(t)\big(E_0(t)\big)$ stays on the same
fixed line $L$, while $h(t)\big(V(t)\big)$ moves along a segment transversally intersecting $L$ at the
point $h(V)$.
Let $U$ be a small (Euclidean)
neighborhood of $V$ in $\Gam$. It follows from Lemma \ref{lreg}, that for each small deformation
of the fragment $(\Gam,\bp,h)\big|_U$ that keeps the image of the deformed edge $E_0$ on the line $L$,
there exists a unique extension up to a deformation of the entire curve $C^*$ such that the images of all marked points but $V$
stay fixed and the combinatorial type of the part $(\Gam,\bp,h)\big|_{\Gam\setminus U}$
does not change. In view of the first formula in (\ref{epsi5}), to prove the constancy of
$RD^{\lab}_y(\bx(t))$, $t\in(t^*-\eps,t^*+\eps)$, it is enough to consider only deformations of
the fragment $(\Gam,\bp,h)\big|_U$. Equivalently, we can assume that
$\Gam$ consists of the unique $(k+3)$-valent vertex $V$ and edges $E_0,...,E_{k+2}$ incident to $V$, while
$\bp$ includes $V$ and one more marked point on the edge $E_0$.

\medskip

{\it (2) The cases $k=0$ and $k=1$.} In this situation, the required statements were proved in \cite[Theorem 1]{IM} and
\cite{GS}, respectively. We provide details for completeness.

If $k=0$, the deformation of $C^*$ in $\widehat{\mathcal M}^{\;e,\lab}_{0,\on}
(\R^2,\Del)$ is presented in Figure \ref{fig2}.
It immediately follows from (\ref{epsi4}) that such a bifurcation does not affect the value of
$RM^{\lab}_y(\bx(t))$, $t\in(t^*-\eps,t^*+\eps)$.

\begin{figure}
\setlength{\unitlength}{1cm}
\begin{picture}(14,3)(0,0)
\thinlines

\put(2,1){\line(1,-1){1}}\put(2,1){\line(1,1){1}}
\put(2.47,1.47){$\bullet$}\put(1.47,0.85){$\bullet$}

\put(12,1){\line(1,-1){1}}\put(12,1){\line(1,1){1}}
\put(12.42,0.38){$\bullet$}\put(11.42,0.85){$\bullet$}

\put(7.5,1){\line(1,-1){0.5}}\put(7.5,1){\line(1,1){0.5}}
\put(7.45,0.87){$\bullet$}\put(6.47,0.85){$\bullet$}

\put(1.9,2){$F_1$}\put(7,2){$F_0$}\put(11.9,2){$F_2$}
\put(4.5,0.9){$\Longleftarrow$}\put(9,0.9){$\Longrightarrow$}

\put(1,1){\line(1,0){1}}\put(6,1){\line(1,0){1.5}}\put(11,1){\line(1,0){1}}

\put(1.3,0.5){$x_1$}\put(6.3,0.5){$x_1$}\put(11.3,0.5){$x_1$}
\put(2.5,1.2){$x_2$}\put(7.8,0.9){$x_2$}\put(12.5,0.7){$x_2$}

\end{picture}
\caption{Degeneration as in Proposition \ref{l1}(2iii) with $k=0$}\label{fig2}
\end{figure}

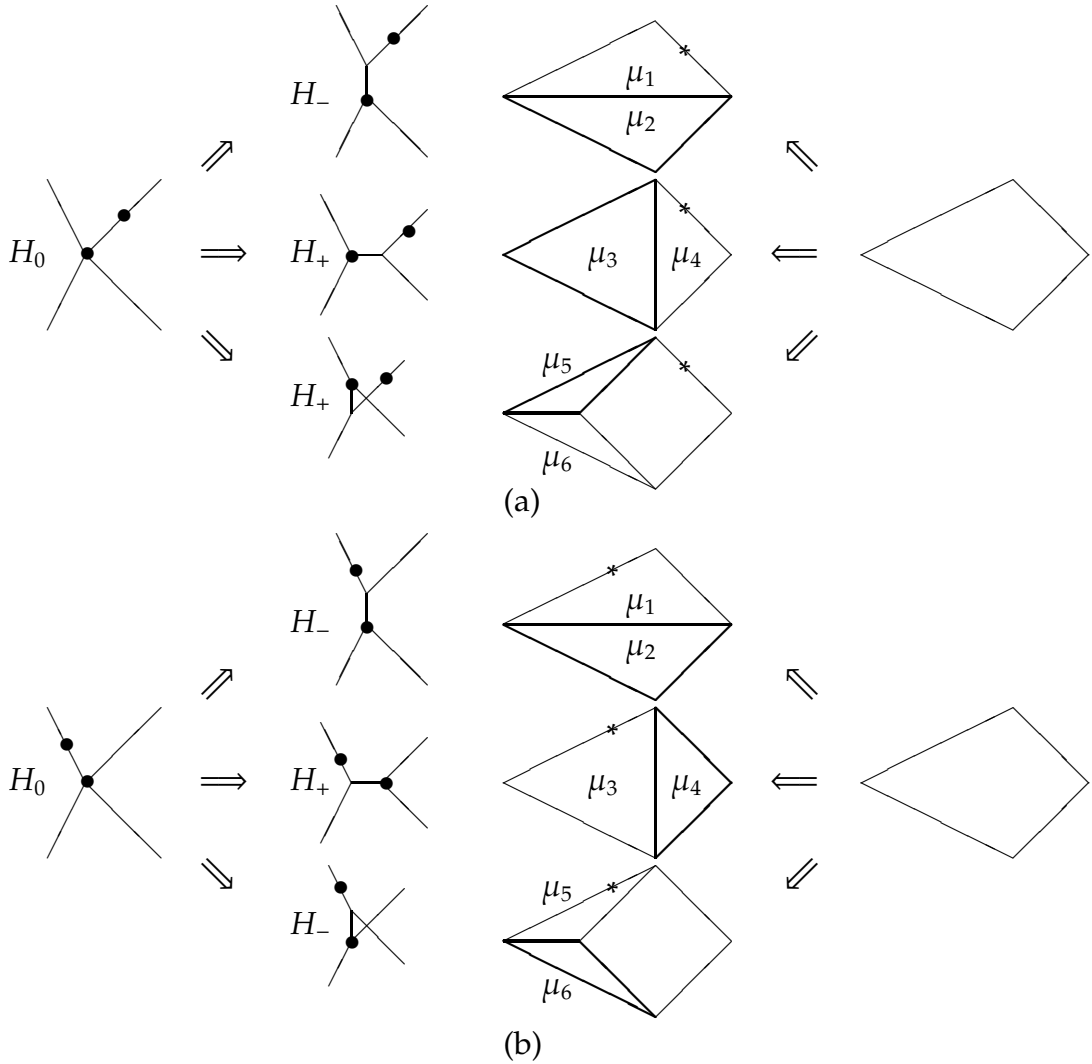
\begin{figure}
\setlength{\unitlength}{1cm}
\begin{picture}(14,14)(0.5,0)
\thinlines

\put(1.5,3.6){\line(1,1){1}}\put(1.5,3.6){\line(1,-1){1}}
\put(1.5,3.6){\line(-1,2){0.5}}\put(1.5,3.6){\line(-1,-2){0.5}}

\put(5.2,5.7){\line(1,-1){0.8}}\put(5.2,5.7){\line(-1,-2){0.4}}
\put(5.2,6.1){\line(1,1){0.8}}\put(5.2,6.1){\line(-1,2){0.4}}
\put(5.2,5.7){\line(0,1){0.4}}

\put(5,3.6){\line(-1,-2){0.4}}\put(5,3.6){\line(-1,2){0.4}}
\put(5.4,3.6){\line(1,1){0.6}}\put(5.4,3.6){\line(1,-1){0.6}}
\put(5,3.6){\line(1,0){0.4}}

\put(5,1.5){\line(-1,-2){0.3}}\put(5,1.9){\line(-1,2){0.3}}
\put(5,1.5){\line(1,1){0.7}}\put(5,1.9){\line(1,-1){0.7}}
\put(5,1.5){\line(0,1){0.4}}

\put(3,3.5){$\Longrightarrow$}
\put(3,4.8){$\Nearrow$}\put(3,2.3){$\Searrow$}

\put(7,0){(b)}\put(7,7.2){(a)}

\put(7,5.7){\line(2,1){2}}
\put(9,6.7){\line(1,-1){1}}

\put(7,3.6){\line(2,1){2}}\put(7,3.6){\line(2,-1){2}}

\put(7,1.5){\line(2,1){2}}
\put(9,2.5){\line(1,-1){1}}\put(9,0.5){\line(1,1){1}}
\put(9,2.5){\line(-1,-1){1}}

\put(11.7,3.6){\line(2,1){2}}\put(11.7,3.6){\line(2,-1){2}}
\put(13.7,4.6){\line(1,-1){1}}\put(13.7,2.6){\line(1,1){1}}

\put(10.5,3.5){$\Longleftarrow$}
\put(10.7,4.8){$\Nwarrow$}\put(10.7,2.3){$\Swarrow$}

\put(0.5,3.5){$H_0$}\put(4.2,5.6){$H_-$}\put(4.2,3.5){$H_+$}\put(4.2,1.6){$H_-$}

\put(8.6,5.9){$\mu_1$}\put(8.6,5.3){$\mu_2$}
\put(8.1,3.5){$\mu_3$}\put(9.2,3.5){$\mu_4$}
\put(7.5,2.1){$\mu_5$}\put(7.5,0.8){$\mu_6$}

\put(1.5,10.6){\line(1,1){1}}\put(1.5,10.6){\line(1,-1){1}}
\put(1.5,10.6){\line(-1,2){0.5}}\put(1.5,10.6){\line(-1,-2){0.5}}

\put(5.2,12.7){\line(1,-1){0.8}}\put(5.2,12.7){\line(-1,-2){0.4}}
\put(5.2,13.1){\line(1,1){0.8}}\put(5.2,13.1){\line(-1,2){0.4}}
\put(5.2,12.7){\line(0,1){0.4}}

\put(5,10.6){\line(-1,-2){0.4}}\put(5,10.6){\line(-1,2){0.4}}
\put(5.4,10.6){\line(1,1){0.6}}\put(5.4,10.6){\line(1,-1){0.6}}
\put(5,10.6){\line(1,0){0.4}}

\put(5,8.5){\line(-1,-2){0.3}}\put(5,8.9){\line(-1,2){0.3}}
\put(5,8.5){\line(1,1){0.7}}\put(5,8.9){\line(1,-1){0.7}}
\put(5,8.5){\line(0,1){0.4}}

\put(3,10.5){$\Longrightarrow$}
\put(3,11.8){$\Nearrow$}\put(3,9.3){$\Searrow$}

\put(7,12.7){\line(2,1){2}}
\put(9,13.7){\line(1,-1){1}}

\put(9,11.6){\line(1,-1){1}}\put(9,9.6){\line(1,1){1}}

\put(7,8.5){\line(2,-1){2}}
\put(9,9.5){\line(1,-1){1}}\put(9,7.5){\line(1,1){1}}
\put(9,7.5){\line(-1,1){1}}

\put(11.7,10.6){\line(2,1){2}}\put(11.7,10.6){\line(2,-1){2}}
\put(13.7,11.6){\line(1,-1){1}}\put(13.7,9.6){\line(1,1){1}}

\put(10.5,10.5){$\Longleftarrow$}
\put(10.7,11.8){$\Nwarrow$}\put(10.7,9.3){$\Swarrow$}

\put(0.5,10.5){$H_0$}\put(4.2,12.6){$H_-$}\put(4.2,10.5){$H_+$}\put(4.2,8.6){$H_+$}

\put(8.6,12.9){$\mu_1$}\put(8.6,12.3){$\mu_2$}
\put(8.1,10.5){$\mu_3$}\put(9.2,10.5){$\mu_4$}
\put(7.5,9.1){$\mu_5$}\put(7.5,7.8){$\mu_6$}

\put(1.42,3.5){$\bullet$}\put(5.1,5.54){$\bullet$}
\put(5.35,3.47){$\bullet$}\put(4.9,1.37){$\bullet$}

\put(1.42,10.5){$\bullet$}\put(5.1,12.54){$\bullet$}
\put(4.9,10.47){$\bullet$}\put(4.9,8.77){$\bullet$}

\put(1.15,4){$\bullet$}\put(4.95,6.3){$\bullet$}\put(4.75,3.8){$\bullet$}\put(4.75,2.1){$\bullet$}
\put(1.9,11){$\bullet$}\put(5.45,13.35){$\bullet$}\put(5.65,10.8){$\bullet$}\put(5.35,8.85){$\bullet$}

\put(9.3,13.2){$*$}\put(9.3,11.1){$*$}\put(9.3,9){$*$}
\put(8.35,6.3){$*$}\put(8.35,4.2){$*$}\put(8.35,2.1){$*$}

\thicklines

\put(7,5.7){\line(2,-1){2}}\put(9,4.7){\line(1,1){1}}
\put(7,5.7){\line(1,0){3}}\put(9,4.6){\line(1,-1){1}}\put(9,2.6){\line(1,1){1}}
\put(9,2.6){\line(0,1){2}}\put(7,1.5){\line(2,-1){2}}
\put(9,0.5){\line(-1,1){1}}
\put(7,1.5){\line(1,0){1}}\put(7,12.7){\line(2,-1){2}}\put(9,11.7){\line(1,1){1}}
\put(7,12.7){\line(1,0){3}}\put(7,10.6){\line(2,1){2}}\put(7,10.6){\line(2,-1){2}}
\put(9,9.6){\line(0,1){2}}\put(7,8.5){\line(2,1){2}}\put(9,9.5){\line(-1,-1){1}}
\put(7,8.5){\line(1,0){1}}

\end{picture}
\caption{Degeneration as in Proposition \ref{l1}(2iii) with $k=1$}\label{fig4}
\end{figure}

Let $k=1$. Suppose that the $h$-images of the four edges incident to $V$ lie of four distinct lines in $\R^2$.
Then $C^*$ admits three types of deformations that correspond to three types of
splitting of the four-valent vertex into a
pair of  trivalent vertices (see Figure \ref{fig4}). We have to study two cases
according as the edge $h(E_0)$
is dual to a side of
the parallelogram inscribed into the quadrangle
or not
(see Figures \ref{fig4}(a,b), where the edge dual to $h(E_0)$ is labeled by asterisk, and the triangles dual to the marked trivalent
vertices are shown
by fat lines).
Thus, (in the notations of Section \ref{sec-plan}) the three top-dimensional cells $F_1,F_2,F_3$ of
$\widehat{\mathcal M}^{\;e,\lab}_{0,\on}(\R^2,\Del)$
attached to $F_0$, project onto $H_+$
or onto $H_-$ according as the moving point $x_2^{(t)}$ belongs to $\R^2_+$ or
$\R^2_-$. In the notation of Section \ref{spsi2}, for the Mikhalkin's multiplicities of the
trivalent vertices
(see Figures \ref{fig4}(a,b)), we have the following additional geometric relations that can be derived from the balancing condition
or by elementary geometry tools
(cf. relations in \cite[Items (B) and (C) in Page 25]{GS})
\begin{equation}\begin{cases}\mu_3=\mu_1+\mu_5,\quad&\text{in Figure \ref{fig4}(a)},\\
\mu_1=\mu_4+\mu_6,\quad&\text{in Figure \ref{fig4}(b)}.\end{cases}\label{e4}\end{equation}
The constancy of
$RD^{\lab}_y(\bx(t))$, $t\in(t^*-\eps,t^*+\eps)$, reduces to the relation
\begin{equation}(z^{\mu_3}-z^{-\mu_3})(z^{\mu_4}+z^{-\mu_4})=
(z^{\mu_1}-z^{-\mu_1})(z^{\mu_2}+z^{-\mu_2})+
(z^{\mu_5}-z^{-\mu_5})(z^{\mu_6}+z^{-\mu_6})\label{epsi6a}\end{equation}
in case of Figure \ref{fig4}(a), and to the relation
\begin{equation}(z^{\mu_1}-z^{-\mu_1})(z^{\mu_2}+z^{-\mu_2})=
(z^{\mu_3}+z^{-\mu_3})(z^{\mu_4}-z^{-\mu_4})+
(z^{\mu_5}+z^{-\mu_5})(z^{\mu_6}-z^{-\mu_6})\label{epsi6b}\end{equation}
in case of Figure \ref{fig4}(b). Both the above equalities
immediately follow from the elementary geometric facts
$$\mu_1+\mu_2=\mu_3+\mu_4,\quad\mu_1-\mu_2=\mu_6-\mu_5,\quad
\mu_3-\mu_4=\mu_5+\mu_6\ ,$$ and formulas (\ref{e4}).

The case when some two edges incident to $V$ are mapped to the same line in $\R^2$
(in Figure \ref{fig4} this means that
the quadrangle turns into a trapeze or a parallelogram), can be treated in the same manner. One can also regard
this case as a limit of the general case considered above and derive the required invariance from the continuous dependence
in the variables $\mu_1,...,\mu_6$
of the expressions in (\ref{epsi6a}) and (\ref{epsi6b}).

\medskip

{\it (3) The case of arbitrary $k\ge1$.}
Let $A=\{\ba_0,...,\ba_{k+2}\}$ be the multiset of the vectors
$$\ba_i:=\ba_V(E_i),\quad i=0,...,k+2\ .$$
Introduce new vectors $$\ba_I=\sum_{i\in I}\ba_i
\quad\text{for any subset}\ I\subset\{0,1,...,k+2\}\ ,$$
and multisets
$$A_I=\{\ba_j,\ 0\le j\le k+2,\ j\not\in I\}\cup\{\ba_I\}
\quad\text{for any subset}\ I\subset\{0,1,...,k+2\}\ ,$$
$$A_{I,J}=\{\ba_j,\ 0\le j\le k+2,\ j\not\in I\cup J\}\cup\{\ba_I,\ba_J\}
\quad\text{for any disjoint subsets}\ I,J\subset\{0,1,...,k+2\}\ .$$
The $C^*$-families bijectively correspond to the vectors $\ba_i$, $1\le i\le k+2$, satisfying $\ba_0\wedge\ba_i\ne0$,
so that, in the deformation
along the family $(\Gam_i(t),\bp_i(t),h_i(t))$, the vertex $V$ splits into an unmarked trivalent vertex
$W_i(t)\in L$, incident to two ends directed by the vectors $\ba_0,\ba_i$, and a marked
$(k+1)$-valent vertex
$V_i(t)$. Since the position of $V_i(t)$ with respect to the line $L$ is determined by
$\sign(\ba_0\wedge\ba_i)$, we can write the required constancy relation in the form
$$\sum_{i=1}^{k+2}
[\ba_0\wedge\ba_i]_y^-\cdot RM_y(\Gam_i(t),\bp_i(t),h_i(t),V_i(t))=0$$
(for any $t\ne t^*$ sufficiently close to $t^*$), or, equivalently (see formula (\ref{epsi2})),
\begin{equation}\sum_{i=1}^{k+2}
[\ba_0\wedge\ba_i]_y^-\cdot \theta_y(A_{\{0,i\}})=0\ .\label{epsi6}\end{equation}
We shall prove relation (\ref{epsi6}) in a greater generality, assuming only that $\sum_{i=0}^{k+2}\ba_i=0$, and we shall use the induction on $k\ge1$.

For $k=1$, relation (\ref{epsi6}) reduces to one of the equalities (\ref{epsi6a}) or (\ref{epsi6b}) proved above. So, assume that $k\ge2$.

Since $|A_{\{i,j\}}|=k+2<k+3=|A|$, the induction assumption yields
\begin{equation}\sum_{\renewcommand{\arraystretch}{0.6}
\begin{array}{c}
\scriptstyle{1\le s\le k+2}\\
\scriptstyle{s\ne i,j}
\end{array}}\left([\ba_0\wedge\ba_s]_y^-
\cdot\theta_y(A_{\{i,j\},\{0,s\}})\right)+[\ba_0\wedge\ba_{\{i,j\}}]_y^-
\cdot \theta_y(A_{\{0,i,j\}})
=0\ .\label{epsi7}\end{equation}
Multiplying (\ref{epsi7}) by $=[\ba_i\wedge\ba_j]_y^+$, we obtain
\begin{equation}\sum_{\renewcommand{\arraystretch}{0.6}
\begin{array}{c}
\scriptstyle{1\le s\le k+2}\\
\scriptstyle{s\ne i,j}
\end{array}}\left([\ba_0\wedge\ba_s]_y^-\cdot [\ba_i\wedge\ba_j]_y^+\cdot\theta_y(A_{\{i,j\},\{0,s\}})\right)+
[\ba_0\wedge\ba_{\{i,j\}}]_y^-\cdot [\ba_i\wedge\ba_j]_y^+\cdot \theta_y(A_{\{0,i,j\}})
=0\ ,\label{epsi8}\end{equation}
Now, for the multiset $\big\{\ba_0,\ba_i,\ba_j,-\ba_{\{0,i,j\}}\big\}$, the induction base yields
$$[\ba_0\wedge\ba_i]_y^-\cdot
[\ba_j\wedge\ba_{\{0,i,j\}}]_y^++
[\ba_0\wedge\ba_j]_y^-\cdot
[\ba_i\wedge\ba_{\{0,i,j\}}]_y^+-
[\ba_0\wedge\ba_{\{0,i,j\}}]_y^-
\cdot[\ba_i\wedge\ba_j]_y^+=0\ ,$$ which after multiplication by
$\theta_y(A_{\{0,i,j\}})$ turns into
$$[\ba_0\wedge\ba_i]_y^-\cdot
[\ba_j\wedge\ba_{\{0,i,j\}}]_y^+\cdot\theta_y(A_{\{0,i,j\}})+
[\ba_0\wedge\ba_j]_y^-\cdot
[\ba_i\wedge\ba_{\{0,i,j\}}]_y^+\cdot\theta_y(A_{\{0,i,j\}})$$
\begin{equation}-
[\ba_0\wedge\ba_{\{0,i,j\}}]_y^-
\cdot[\ba_i\wedge\ba_j]_y^+\cdot\theta_y(A_{\{0,i,j\}})=0\ .\label{epsi10}\end{equation}
We sum up relations (\ref{epsi8}) and (\ref{epsi10}), observing that the last summands in their left-hand side cancel out in view of $\ba_0\wedge\ba_{\{0,i,j\}}=\ba_0\wedge\ba_{\{i,j\}}$. Then we sum up the resulting equalities for all pairs $(i,j)$
such that $1\le i<j\le k+2$ and obtain the following:
$$0=\sum_{1\le i<j\le k+2}\sum_{\renewcommand{\arraystretch}{0.6}
\begin{array}{c}
\scriptstyle{1\le s\le k+2}\\
\scriptstyle{s\ne i,j}
\end{array}}[\ba_0\wedge\ba_s]_y^-\cdot [\ba_i\wedge\ba_j]_y^+\cdot\theta_y(A_{\{i,j\},\{0,s\}})$$
$$+\sum_{1\le i<j\le k+2}\left([\ba_0\wedge\ba_i]^-_y\cdot[\ba_j\wedge\ba_{\{0,i,j\}}]^+_y+
[\ba_0\wedge\ba_j]^-_y\cdot[\ba_i\wedge\ba_{\{0,i,j\}}]^+_y\right)\theta_y(A_{\{0,i,j\}})$$
$$=\sum_{s=1}^{k+2}\Bigg([\ba_0\wedge\ba_s]^-_y\cdot\sum_{\renewcommand{\arraystretch}{0.6}
\begin{array}{c}
\scriptstyle{1\le i<j\le k+2}\\
\scriptstyle{i,j\ne s}\end{array}}[\ba_i\wedge\ba_j]^+_y\cdot\theta_y(A_{\{0,i,j\}})\Bigg)$$
$$+\sum_{s=1}^{k+2}\bigg([\ba_0\wedge\ba_s]^-_y\cdot\sum_{s<j\le k+2}
[\ba_j\wedge\ba_{\{0,s\}}]^+_y\cdot\theta_y(A_{\{0,s,j\}})\Bigg)$$
$$+\sum_{s=1}^{k+2}\bigg([\ba_0\wedge\ba_s]^-_y\cdot\sum_{1\le i<s}
[\ba_i\wedge\ba_{\{0,s\}}]^+_y\cdot\theta_y(A_{\{0,s,i\}})\Bigg)$$
$$=\sum_{s=1}^{k+2}[\ba_0\wedge\ba_s]_y^-\cdot
\Bigg(\sum_{\renewcommand{\arraystretch}{0.6}
\begin{array}{c}
\scriptstyle{1\le i<j\le k+2}\\
\scriptstyle{i,j\ne s}
\end{array}}[\ba_i\wedge\ba_j]_y^+\cdot \theta_y(A_{\{i,j\},\{0,s\}})
+\sum_{\renewcommand{\arraystretch}{0.6}
\begin{array}{c}
\scriptstyle{1\le i\le k+2}\\
\scriptstyle{i\ne s}
\end{array}}[\ba_i\wedge\ba_{\{0,s\}}]_y^+\cdot \theta_y(A_{\{0,s,i\}})
\Bigg)\ .$$
It remains to notice that by definition (\ref{emu}), the expression in the parentheses
equals $\theta_y(A_{\{0,s\}})$. \proofend

\subsection{On the uniqueness of the refinement}\label{secu}
The summands in the right-hand side of the recursion (\ref{emu}) are enumerated by splittings of
a tree with a unique $(k+2)$-valent vertex into trees with a $(k+1)$-valent vertex and a trivalent vertex.
One, however, can use a recursion based on splittings into trees of other type
(for instance, into trees having two trivalent vertices and a $k$-valent vertex etc.).
Iterating such a recursion, we finally end up with the sum
of the form
$\sum_\alpha \left(c_\alpha\sum_T\mu_y(T)\right)$, where $\alpha$ runs over the set of combinatorial types of trivalent trees with $k+2$ leaves, $T$ ranges over all possible labelings of the leaves of a fixed tree of type $\alpha$ by elements of the sequence $A$,
and $\mu_y(T)$ is the product of the factors $[\mu(V)]^+_y$ over all the vertices of $T$.

We intend to show that any
such definition leads to a refined invariant that differs from $RD_y(\Del,\on)$ by a
constant multiplicative factor depending only on combinatorics of splittings used in the recursion. In this sense we speak of the uniqueness of the refinement of rational descendant tropical invariants.

Let $\Del$ be a balanced sequence of $m\ge3$ vectors in $\R^2$. Let $\alpha$ be a trivalent tree with $m$ leaves,
$\Leaf(\alpha)$ the set of its leaves, $\alpha^0$ the set of its trivalent vertices.
Let $\varphi:\Leaf(\alpha)\to\Del$ be a bijection.
It is easy to show that there exists a unique map $\psi_\varphi:\Fl(\alpha)\to\R^2$, where $\Fl(\alpha)$ is the set of pairs
$(V,E)$ with $V\in\alpha^0$ and $E$ an edge incident to $V$, such that
\begin{itemize}\item if $E\in\Leaf(\alpha)$ then $\psi_\varphi(V,E)=\varphi(E)$,
\item if $E$ is incident to two vertices $V_1,V_2\in\alpha^0$, then
$\psi_\varphi(V_1,E)+\psi_\varphi(V_2,E)=0$,
\item for any vertex $V\in\alpha^0$ and $E_1,E_2,E_3$ edges of $\alpha$ incident to $V$,
$$\psi_\varphi(V,E_1)+\psi_\varphi(V,E_2)+\psi_\varphi(V,E_3)=0\ .$$\end{itemize}
Denote the triple of vectors in the latter relation by $\Del_{\alpha,\varphi}(V)$. Define
\begin{equation}\theta_{y,\alpha}(\Del)=\sum_{\varphi}\prod_{V\in\alpha^0}
\theta_y(\Del_{\alpha,\varphi}(V))\ ,\label{euni1}\end{equation}
where $\varphi$ ranges over all bijections $\Leaf(\alpha)\to\Del$.

\begin{lemma}\label{l-uni}
For any $m\ge3$, any balanced sequence $\Del$ of $m$ vectors in $\R^2$, and any trivalent trees $\alpha,\beta$
with $m$ leaves, one has
$\theta_{y,\alpha}(\Del)=\theta_{y,\beta}(\Del)$.
\end{lemma}

{\bf Proof.}
Introduce the following notation: for any sequence of vectors $\bb_1,...,\bb_r$ ($r\ge2$) and a permutation
$\sigma\in S_r$, put $\sigma\Lambda(\bb_1,...,\bb_r)=
\sum_{1\le i<j\le r}\bb_{\sigma(i)}\wedge\bb_{\sigma(j)}$.

Given a trivalent tree $\alpha$ with $m$ leaves, for any $\ba\in\Del$ and any $E\in\Leaf(\alpha)$, define
\begin{equation}\theta_{y,\alpha}(\Del,\ba,E)=\sum_{\varphi(E)=\ba}\prod_{V\in\alpha^0}
\theta_y(\Del_{\alpha,\varphi}(V))\ ,\label{ne23}\end{equation}
where $\varphi$ ranges over all bijections
$\Leaf(\alpha)\to\Del$ satisfying $\varphi(E)=\ba$.
In view of the relations
\begin{equation}
\theta_{y,\alpha}(\Del)=\sum_{\ba\in\Del}\theta_{y,\alpha}(\Del,\ba,E)\quad\text{for all}\ E\in\Leaf(\alpha)\ ,
\label{e3}\end{equation} the following claim completes the proof of Lemma:
For any $\ba\in\Del$ and $E\in\Leaf(\alpha)$, the following holds
\begin{equation}\theta_{z^2,\alpha}(\Del,\ba,E)=\frac{2^{m-2}}{(z+z^{-1})^{m-2}}
\sum_{\sigma\in S_{m-1}}z^{\sigma\Lambda(\Del\setminus\{\ba\})}\ ,
\label{e2}\end{equation}
where the right-hand side does not depend neither on the choice of $E$ nor on the combinatorial type of $\alpha$.

We prove formula (\ref{e2}) by induction on $m$. For $m=3$, formula (\ref{e2}) immediately follows from (\ref{ne23}).
Suppose that $m\ge4$.

If $E$ is incident to a trivalent vertex together with another edge $E'\in\Leaf(\alpha)$,
then by definition (\ref{ne23}),
we obtain
$$\theta_{z^2,\alpha}(\Del,\ba,E)=\sum_{\bb\in\Del\setminus\{\ba\}}
\frac{z^{\ba\wedge\bb}+z^{\bb\wedge\ba}}{z+z^{-1}}
\cdot\theta_{z^2,\gamma}(\widetilde\Del,\ba+\bb,\widetilde E)\ ,$$
where the trivalent tree $\gamma$ with $m-1$ leaves is obtained from $\alpha$ by removing the leaves $E,E'$,
the sequence $\widetilde\Del$ is obtained from $\Del$ by removing $\ba$ and $\bb$ and adding the vector $\ba+\bb$ (say, as the last vector), and
$\widetilde E\in\Leaf(\gamma)$ is the third edge of $\alpha$ incident to $V$. Using the induction assumption, we obtain
$$\theta_{z^2,\alpha}(\Del,\ba,E)=\sum_{\bb\in\Del\setminus\{\ba\}}
\frac{z^{\ba\wedge\bb}+z^{\bb\wedge\ba}}{z+z^{-1}}
\cdot\frac{2^{m-3}}{(z+z^{-1})^{m-3}}\sum_{\sigma\in S_{m-2}}z^{\sigma\Lambda(\Del\setminus\{\ba,\bb\})}$$
\begin{equation}=\frac{2^{m-3}}{(z+z^{-1})^{m-2}}\sum_{\bb\in\Del\setminus\{\ba\}}
\sum_{\sigma\in S_{m-2}}\left(z^{\ba\wedge\bb+\sigma\Lambda(\Del\setminus\{\ba,\bb\}}
+z^{\bb\wedge\ba+\sigma\Lambda(\Del\setminus\{\ba,\bb\}}\right)\ .\label{e2a}\end{equation}
Since $\ba=-\bb-\sum_{\bc\in\Del\setminus\{\ba,\bb\}}\bc$, the exponents of $z$ in the latter expression turn into
$$\ba\wedge\bb+\sigma\Lambda(\Del\setminus\{\ba,\bb\}=\sum_{\bc\in\Del\setminus\{\ba,\bb\}}\bb\wedge\bc
+\sigma\Lambda(\Del\setminus\{\ba,\bb\})=
\sigma^{\bb}\Lambda(\Del\setminus\{\ba\})\ ,$$
$$\bb\wedge\ba+\sigma\Lambda(\Del\setminus\{\ba,\bb\}=\sum_{\bc\in\Del\setminus\{\ba,\bb\}}\bc\wedge\ba
+\sigma\Lambda(\Del\setminus\{\ba,\bb\})=
\sigma_{\bb}\Lambda(\Del\setminus\{\ba\})\ ,$$ where permutations $\sigma^{\bb},\sigma_{\bb}$ of $\Del\setminus\{\ba\}$ are obtained from $\sigma$ by sending $\bb$ to the first or to the last place, respectively. It follows that, when $\bb$ ranges over all vectors of $\Del\setminus\{\ba\}$ and $\sigma$ ranges over all permutations of
$\Del\setminus\{\ba,\bb\}$, the permutations $\sigma^{\bb}$ and $\sigma_{\bb}$ twice run over all permutations of $\Del\setminus\{\ba\}$, and hence (\ref{e2a}) yields the required relation (\ref{e2}).

Let $E$ be the only leaf incident to a vertex $V\in\alpha^0$. Removing $E$ from $\alpha$, we obtain two trees
$\alpha_1$, $\alpha_2$ sharing the vertex $V$, with $m_1$ and $m_2$ leaves, respectively, where $m_1,m_2\ge3$ and $m_1+m_2=m+1$. Denote by
$E_1$, $E_2$ the leaves of $\alpha_1$, $\alpha_2$, respectively, incident to $V$. Hence
$$\theta_{z^2,\alpha}(\Del,\ba,E)=\sum_{\renewcommand{\arraystretch}{0.6}
\begin{array}{c}
\scriptstyle{I_1\cup I_2=\Del\setminus\{\ba\}}\\
\scriptstyle{I_1\cap I_2=\emptyset}\\
\scriptstyle{|I_1|=m_1-1,\ |I_2|=m_2-1}
\end{array}}\frac{z^{\bb_1\wedge\bb_2}+z^{\bb_2\wedge\bb_1}}{z+z^{-1}}
\cdot\theta_{y,\alpha_1}(\Del_1,\bb_1,E_1)
\cdot\theta_{y,\alpha_2}(\Del_2,\bb_2,E_2)\ ,$$ where
$$\bb_i=-\sum_{\bb\in I_i}\bb,\quad\Del_i=I_i\cup\{\bb_i\},\quad i=1,2\ .$$
Since $m_1,m_2<m$, we apply the induction assumption and obtain
$$\theta_{z^2,\alpha}(\Del,\ba,E)=\sum_{\renewcommand{\arraystretch}{0.6}
\begin{array}{c}
\scriptstyle{I_1\cup I_2=\Del\setminus\{\ba\}}\\
\scriptstyle{I_1\cap I_2=\emptyset}\\
\scriptstyle{|I_1|=m_1-1,\ |I_2|=m_2-1}
\end{array}}\frac{z^{\bb_1\wedge\bb_2}+z^{\bb_2\wedge\bb_1}}{z+z^{-1}}
\cdot\frac{2^{m-3}}{(z+z^{-1})^{m-3}}\cdot\sum_{\sigma_1\in S_{m_1-1}}z^{\sigma_1\Lambda(I_1)}\cdot\sum_{\sigma_2\in S_{m_2-1}}
z^{\sigma_2\Lambda(I_2)}$$
$$=\frac{2^{m-3}}{(z+z^{-1})^{m-2}}\sum_{\renewcommand{\arraystretch}{0.6}
\begin{array}{c}
\scriptstyle{I_1\cup I_2=\Del\setminus\{\ba\}}\\
\scriptstyle{I_1\cap I_2=\emptyset}\\
\scriptstyle{|I_1|=m_1-1,\ |I_2|=m_2-1}
\end{array}}\sum_{\renewcommand{\arraystretch}{0.6}
\begin{array}{c}
\scriptstyle{\sigma_1\in S_{m_1-1}}\\
\scriptstyle{\sigma_2\in S_{m_2-1}}
\end{array}}\left(z^{\bb_1\wedge\bb_2+\sigma_1\Lambda(I_1)+\sigma_2\Lambda(I_2)}
+z^{\bb_2\wedge\bb_1+\sigma_1\Lambda(I_1)+\sigma_2\Lambda(I_2)}\right)\ .$$
Since $\bb_1=-\sum_{\bb\in I_1}\bb$ and $\bb_2=-\sum_{\bc\in I_2}\bc$, we get
$$\bb_1\wedge\bb_2+\sigma_1\Lambda(I_1)+\sigma_2\Lambda(I_2)=\sigma_{12}\Lambda(\Del\setminus\{\ba\})\ ,$$ $$\bb_2\wedge\bb_1+\sigma_1\Lambda(I_1)+\sigma_2\Lambda(I_2)=\sigma_{21}\Lambda(\Del\setminus\{\ba\})\ ,$$ where $\sigma_{12}$, resp. $\sigma_{21}$, is a permutation of $\Del\setminus\{\ba\}$ obtained from $\sigma_1,\sigma_2$ by setting the elements of $I_1$ before $I_2$, resp. after $I_2$. It follows that each of the permutations $\sigma_{12}$ and $\sigma_{21}$ ranges over the whole group of permutations of
$\Del\setminus\{\ba\}$, which finally yields formula (\ref{e2}).
\proofend

\section{Examples}\label{sec-ex}

In this section we calculate the invariant $RM_y^{lab}(\Del)$
in a series of examples, in which the upper bound (\ref{ed1}) to the degree of
the denominator is sharp.

The first example is degenerate and plays an auxiliary role.

\begin{lemma}\label{al1}
Let $m_1,...,m_r$ be positive integers, where $r\ge2$, and $m=m_1+...+m_r$. Set
$\Del^0_r=\{(m_1,0),\ ...,\ (m_r,0),\ (-a,0)\}$ and $\on_r=(n_k)_{k\ge0}$ such that $n_{r-1}=1$,
$n_k=0$, $k\ne r-1$ (i.e., consider tropical curves with a marked $(r+1)$-valent vertex). Then
$$RD_y^{\lab}(\Del^0_r,\on_r)=\frac{r!\cdot (r+1)!}{6(y^{1/2}+y^{-1/2})^{r-1}}\ .$$
\end{lemma}

{\bf Proof.} The formula evidently holds for $r=2$, and it can be proved by induction on $r$ using the recursive formula
$$RD_y^{\lab}(\Del^0_r,\on_r)=\binom{r+1}{2}\cdot\frac{2}{y^{1/2}+y^{-1/2}}\cdot RD^{\lab}_y(\Del^0_{r-1},\on_{r-1})$$
that immediately follows from (\ref{emu}) and (\ref{epsi2}).
\proofend

The next example deals with a nondegenerate Newton triangle.

\begin{figure}[ht]	
	\centering
	\begin{tikzpicture}	
        \draw (0,1) -- (2,0) -- (2,5) --  cycle;
        \draw node at (0,1) [regular_node]{};
        \draw node at (2,0) [regular_node]{};
        \draw node at (2,5) [regular_node]{};
        \draw node at (2,1) [small_node]{};
        \draw node at (2,2.25) [small_node]{};
        \draw node at (2,3.5) [small_node]{};
        \draw node at (2,4) [small_node]{};
        \node at (2.5,0.5) {$m_1$};
        \node at (2.5,1.7) {$m_2$};
        \node at (2.5,2.7) {$\vdots$};
        \node at (2.5,4.5) {$m_r$};
        \node at (2.5,3.75) {$m_{r-1}$};
    \end{tikzpicture}
     \caption{Newton polygon $P(\Del^h_{m_1...m_r})$}\label{fig8}
\end{figure}
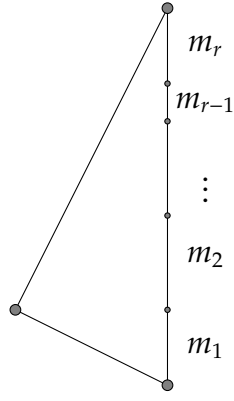

\begin{lemma}\label{al2}
Let $h$ and $m_1,...,m_r$ be positive integers, where $r\ge2$, and $m=m_1+...,m_r$. Set
$$\Del^h_{m_1...m_r}=\{(m_1,0),\ ...,\ (m_r,0),\ (-1,-h),\ (1-m,h)\}\quad\text{(cf. Figure
\ref{fig8})}\ .$$
Then
$$RD_y^{lab}(\Del^h_{m_1\ldots m_r},\on_r)=\frac{(r+2)!}{12\cdot (y^{1/2}+y^{-1/2})^r}
\qquad\qquad\qquad$$
\begin{equation}\qquad\qquad\qquad\times\sum_{I}(r-|I|)!\cdot|I|!\cdot\left\{y^{h(m/2-\sum_{i\in I}{m_i})}+y^{h(\sum_{i\in I}{m_i}-m/2)}\right\}\ ,\label{ae1}\end{equation}
where the sum runs over all subsets $I\subset\{1,2,...,r\}$.
\end{lemma}

\begin{figure}[h]	
	\centering
	\begin{tikzpicture}	
        \draw (0,1) -- (2,0) -- (2,5) --  cycle;
        \draw node at (0,1) [regular_node]{};
        \draw node at (2,0) [regular_node]{};
        \draw node at (2,5) [regular_node]{};
        \node at (1.5,2.5) {$\sum {m_i}$};
    \end{tikzpicture}
    \begin{tikzpicture}
        \draw (1.9,0) -- (1.9,5) --(2,5) --(2,0) -- cycle;
        \draw node at (2,0) [regular_node]{};
        \draw node at (2,5) [regular_node]{};
        \draw node at (2,1) [small_node]{};
        \draw node at (2,2.25) [small_node]{};
        \draw node at (2,3.5) [small_node]{};
        \draw node at (2,4) [small_node]{};
        \node at (2.5,0.5) {$m_1$};
        \node at (2.5,1.7) {$m_2$};
        \node at (2.5,2.7) {$\vdots$};
        \node at (2.5,4.5) {$m_r$};
        \node at (2.5,3.75) {$m_{r-1}$};
    \end{tikzpicture} \hspace{2cm}
    \begin{tikzpicture}	
        \draw (0,1) -- (2,5) -- (2,4) --  cycle;
        \draw node at (0,1) [regular_node]{};
        \draw node at (2,5) [regular_node]{};
        \draw node at (2,4) [regular_node]{};
        \node at (2.3,4.5) {$ {m_i}$};
        \draw (0.3,1) -- (2.3,0) -- (2.3,4) --  cycle;
        \draw node at (2.3,0) [regular_node]{};
        \draw node at (0.3,1) [regular_node]{};
        \draw node at (2.3,0.8) [small_node]{};
        \draw node at (2.3,1.75) [small_node]{};
        \draw node at (2.3,2.15) [small_node]{};
        \draw node at (2.3,2.75) [small_node]{};
        \draw node at (2.3,3.5) [small_node]{};
        \draw node at (2.3,4) [small_node]{};
        \node at (2.8,0.5) {$m_1$};
        \node at (2.8,1.2) {$m_2$};
        \node at (2.8,2) {$\vdots$};
        \node at (2.8,2.5) {$\widehat{m_i}$};
        \node at (2.8,3.25) {$\vdots$};
        \node at (2.8,3.75) {$m_r$};
    \end{tikzpicture} \hspace{2cm}
    \begin{tikzpicture}
        \draw (0,1) -- (2,0) -- (2,5) --  cycle;
        \draw (2.3,4)--(2.4,4) -- (2.4,5) --(2.3,5) -- cycle;
        \draw node at (0,1) [regular_node]{};
        \draw node at (2,0) [regular_node]{};
        \draw node at (2,5) [regular_node]{};
        \draw node at (2,0.8) [small_node]{};
        \draw node at (2,1.15) [small_node]{};
        \draw node at (2,1.75) [small_node]{};
        \draw node at (2,2.15) [small_node]{};
        \draw node at (2,2.75) [small_node]{};
        \draw node at (2,3.5) [small_node]{};
        \draw node at (2,4) [small_node]{};
        \draw node at (2.4,4) [regular_node]{};
        \draw node at (2.4,5) [regular_node]{};
        \draw node at (2.4,4.5) [small_node]{};
        \node at (1.5,0.5) {$m_1$};
        \node at (1.5,1) {$\vdots$};
        \node at (1.5,1.5) {$\widehat{m_i}$};
        \node at (1.5,2.15) {$\vdots$};
        \node at (1.5,2.55) {$\widehat{m_j}$};
        \node at (1.5,3.2){$\vdots$};
        \node at (1.5,3.75) {$m_r$};
        \node at (1.1,4.55){$m_i+m_j$};
        \node at (2.8,4.25) {$m_i$};
        \node at (2.8,4.75){$m_j$};
    \end{tikzpicture}
     \caption{Splittings in formula (\ref{e19})}\label{fig9}
\end{figure}
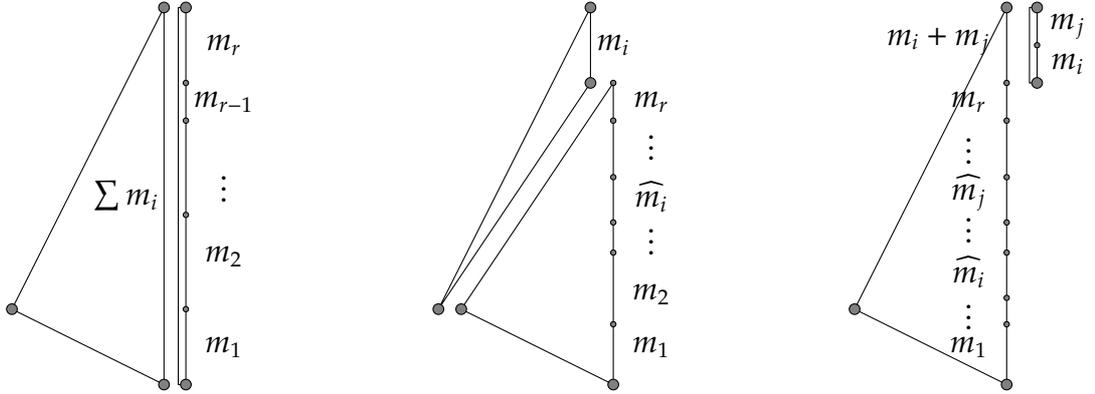

{\bf Proof.} Formula (\ref{epsi2}) yields the relation
\begin{equation}
    \begin{aligned}
       RD_y^{lab}(\Del^h_{m_1\ldots m_r},\on_r)={} & \left\{y^{mh/2}+y^{-mh/2}\right\}\cdot RD_y^{lab}(\Del^0_r,\on_r) +\\
         &+2\sum_{i=1}^r \frac{y^{m_ih/2}+y^{-m_ih/2}}{y^{1/2}+y^{-1/2}}\cdot RD_y^{lab}(\Del^h_{m_1\ldots\widehat{m_i}\ldots m_r},\on_{r-1})+\\
         &+\sum_{i\neq j}\frac {2}{y^{1/2}+y^{-1/2}}\cdot RD_y^{lab}(\Del^h_{m_1\ldots\widehat{m_i}\ldots\widehat{m_j}\ldots m_r,m_i+m_j},
         \be_{r-1})
    \end{aligned}\label{e19}
\end{equation} that corresponds to the three types of splittings shown in Figure \ref{fig9}.

A routine induction on $r$, skipped here for brevity,
completes the proof.
\proofend

\begin{corollary}\label{ac1}
(1) Let, in the notation of Lemma \ref{al2}, $r=2k\ge2$, $m_1=...=m_{2k}=1$, and $h=2$. Then
\begin{equation}RD_y^{\lab}(\R^2,\Del^2_{2k\times1})=\frac{(2k)!\cdot(2k+2)!}{6}\cdot\frac{y^{2k}
+y^{2k-2}+...+y^{2-2k}+y^{-2k}}
{(y+2+y^{-1})^k}\ .\label{ec3a}\end{equation}

(2) Let, in the notation of Lemma \ref{al2}, $r=2k+1\ge3$, $m_1=...=m_{2k+1}=1$, and $h>1$ satisfy $\gcd(h,2k)=1$. Then
\begin{equation}RD_y^{\lab}(\R^2,\Del^h_{(2k+1)\times1})
=\frac{(2k+3)!\cdot(2k+1)!}{6(y+2+y^{-1})^k}\cdot
\sum_{i=0}^{h-1}(-1)^iy^{(h-1)/2-i}\cdot\sum_{i=-k}^ky^{hi}\ .\label{ec3b}\end{equation}
\end{corollary}

This is a consequence of formula (\ref{ae1}). Note that the bound (\ref{ed1}) to the degree of denominators of refined descendant invariants turns into an equality
under the hypotheses of Corollary \ref{ac1}, since the denominator in formulas (\ref{ec3a}) and (\ref{ec3b}) is coprime to other terms.

\section*{Appendix: Lattice path algorithm}\label{sec-lp}

In \cite{Mi} Mikhalkin proved that the Gromov-Witten invariants of toric surfaces
can be computed by summing up multiplicities of finitely many specific trivalent plane tropical curves
with marked points on edges, and he suggested a combinatorial algorithm (so-called lattice path algorithm), which associated the counted tropical curves
with certain subdivisions of the given Newton polygon into convex lattice polygons, while the Mikhalkin's multiplicity
of a tropical curve in count appeared to be the product of multiplicities of the pieces of the corresponding subdivision.
H. Markwig and J. Rau \cite[Section 9]{MR} generalized this algorithm to computation of descendant rational tropical Gromov-Witten invariants:
now replacing rational plane tropical curves with marked points on edges and at vertices by suitable
subdivisions of the Newton polygon and computing the multiplicity of each curve as the product of multiplicities of the elements of the corresponding subdivision. The same algorithm applies well for the computation of our refined descendant invariant, since we enumerate the same collection of
rational marked tropical curves and our refined multiplicity is also the product of multiplicities of the vertices.

For the background of the algorithm we refer to \cite[Section 9]{MR}.
Here we just describe it
and explain how to compute the refined descendant invariant.
As illustration, we consider the example studied in \cite[Example 9.28]{MR}.

Let $\Del\subset\Z^2\setminus\{0\}$ be a nondegenerate balanced multi-set.
Let $\lambda:\R^2\to\R$ be a linear functional injective on $\Z^2$. Orient each line $\lambda=\const$ so that after the clockwise rotation by $\frac{\pi}{2}$ it becomes $\lambda$-ordered.
Fix the linear $\lambda$-order on the set $P(\Del\cap\Z^2)$
denoted by $\prec$.
Denote by $q_{\min}$ and $q_{\max}$ the minimal and the maximal points in $P(\Del)$. A lattice $\lambda$-path of length $m$ in
$P(\Del)$ is a broken line with vertices $q_0\prec q_1\prec...\prec q_m$ such that $q_0=q_{\min}$, $q_m=q_{\max}$.
Each lattice $\lambda$-path $\gamma$ divides the strip $\Sigma=\{\lambda(q_{\min})\le\lambda(x,y)\le\lambda(q_{\max})\}$ into two parts,
whose closures we denote $\Sigma^+(\gamma)$ and $\Sigma^-(\gamma)$ in accordance with the orientation of the lines $\lambda=\const$. The algorithm consists of
the three procedures:
\begin{enumerate}\item[(A)] construction of the initial data, two lattice $\lambda$-paths, $\gamma^+_0,\gamma^-_0$ in
$P(\Del)$ such that $\gamma^-_0\subset\Sigma^-(\gamma^+_0)$ and $\gamma^+_0\subset\Sigma^+(\gamma^-_0)$;
\item[(B)] subdivision of the domains $\Sigma^+(\gamma^+_0)\cap P(\Del)$ and $\Sigma^-(\gamma^-_0)\cap P(\Del)$ into
lattice triangles and parallelograms;
\item[(C)] computation of the refined weight of each admissible
subdivision resulting from the two preceding procedures and summing up over all these
subdivisions.
\end{enumerate}

For the procedure (A), take any sequence $m_1,m_2,...,m_n\ge2$, in which any value $m\ge2$ is attained
exactly $n_{m-2}$ times, and take arbitrary partitions
$$m_i=m'_i+m''_i\quad\text{such that}\quad m'_i,m''_i\ge1,\quad i=1,...,n\ .$$ Then take two lattice $\lambda$-paths $\gamma^+_0$,
having $\sum_im'_i+1$ vertices, and $\gamma^-_0$, having $\sum_im''_i+1$ vertices, and such that
\begin{itemize}\item $\gamma^-_0\subset\Sigma^-(\gamma^+_0)$, $\gamma^+_0\subset\Sigma^+(\gamma^-_0)$,
\item for each $k=1,...,n-1$, the vertex $v^+_{s'(k)}$ of $\gamma^+_0$ coincides with the vertex
$v^-_{s''(k)}$ of $\gamma^-_0$, where $s'(k)=\sum_{i\le k}m'_i$ and $s''(k)=\sum_{i\le k}m''_i$,
\item for each $k=1,...,n$, the convex hull $Q_k$ (called a {\it rag rug element} in \cite{MR}) of the
vertices $v^+_i$, $s'(k-1)\le i\le s'(k)$, of $\gamma^+_0$ together
with the vertices $v^-_i$, $s''(k-1)\le i\le s''(k)$, of $\gamma^-_0$ is such that all the aforementioned
vertices belong to $\partial Q_k$.
\end{itemize}

The procedure (B) is the same as in \cite[Section 7.2]{Mi}. It produces two sequences of lattice $\lambda$-paths,
$\gamma^+_k$, $k\ge0$, and $\gamma^-_k$, $k\ge0$. Given a lattice $\lambda$-path $\gamma^+_k$, resp.
$\gamma^-_k$, with vertices $q_0\prec q_1\prec...\prec q_m$, we
take
$$j=\min\left\{1\le i<m,\ [q_{i-1},q_{i+1}]\subset\Sigma^+(\gamma^+_k),\ q_i\not\in[q_{i-1},q_{i+1}]
\right\}\ ,$$
resp.
$$j=\min\left\{1\le i<m,\ [q_{i-1},q_{i+1}]\subset\Sigma^-(\gamma^-_k),\ q_i\not\in[q_{i-1},q_{i+1}]
\right\}\ .$$
If such $j$ does not exist, we say that the path $\gamma^+_k$, resp.
$\gamma^-_k$, is terminal. If such $j$ exists, we define the lattice $\lambda$-path
$\gamma_{k+1}$ (resp. $\gamma^-_{k+1}$)
\begin{itemize}\item either by the sequence of vertices $q_0,...,q_{j-1},q_{j+1},...,q_m$; in this case we
include the triangle $\conv(q_{j-1},q_j,q_{j+1})$ into the set of tiles of the constructed subdivision,
\item or by the sequence $q_0,...,q_{j-1},q'_j,q_{j+1},...,q_m$, where $\conv(q_{j-1},q_jq'_j,q_{j+1})$ is a parallelogram,
provided that $q'_j\in P(\Del)$; in this case we include the parallelogram
$\conv(q_{j-1},q_jq'_j,q_{j+1})$ into the set of tiles
of the constructed subdivision of $P(\Del)$. \end{itemize} Let $\gamma^+_k$ and $\gamma^-_l$ be terminal.
We say that the obtained subdivision of $P(\Del)$ is {\it admissible} if $\gamma^+_k\cup\gamma^-_l=\partial P(\Del)$,
and it meets
\begin{itemize}\item the degree condition: the segments of $\gamma^+_k$ and $\gamma^-_l$, being rotated by $\frac{\pi}{2}$ and
oriented as exterior normal vectors to $\partial P(\Del)$, form the entire
multiset $\Del$, and
\item the connectedness condition of the following dual graph $G$: each rag rug element $Q_i$, $1\le i\le n$, and
each triangle in $\Sigma^+(\gamma^+_0)\cup\Sigma^-(\gamma^-_0)$ corresponds to a vertex of $G$, each parallelogram
corresponds to two vertices that are associated with two pairs of parallel sides; two vertices,
coming from polygons in $\Sigma^+(\gamma^+_0)$ (resp., $\Sigma^-(\gamma^-_0)$) are joined by an arc
if the corresponding polygons share a common side, a vertex corresponding to $Q_i$ and a vertex corresponding to a polygon
in $\Sigma^+(\gamma^+_0)$ (resp., $\Sigma^-(\gamma^-_0)$) are joined by an arc if $\gamma^+_0$
(resp., $\gamma^-_0$) and the polygon share a segment.\end{itemize}

The refined multiplicity of an admissible subdivision is of the form $Z_1Z_2/Z_3$. Here, $Z_1$ is the product of the factors
$[\nu]^-_y$ over all triangles in $\Sigma^+(\gamma^+_0)\cup\Sigma^-(\gamma^-_0)$ that are $[\nu]^-_y$, where $\nu$ is the lattice area of the triangle.
The term $Z_2$ is the product of the refined multiplicities of all rag rug elements:
if $$Q_k=\conv(\{v'_i\}_{s'(k-1)\le i\le s'(k)},\{v''_i\}_{s''(k-1)\le i\le s''(k)})\ ,$$ we build a multiset of vectors
$\Del_k$ obtained from the segments $[v'_i,v'_{i+1}]$, $s'(k-1)\le i<s'(k)$, rotated by $\frac{\pi}{2}$ and oriented towards $\Sigma^+(
\gamma^+_0)$, and from the segments $[v''_i,v''_{i+1}]$, $s''(k-1)\le i<s''(k)$, rotated by $\frac{\pi}{2}$ and oriented towards
$\Sigma^-(\gamma^-_0)$, and we set the refined multiplicity of $Q_k$ to be $\theta_y(\Del_k)$.
At last, $Z_3$ is the product of orders of the automorphisms groups of $\Del_1,....,\Del_n$.
The invariant $RD_y(\Del,\on)$ equals the sum of the refined multiplicities of all the admissible subdivisions resulting from the
above algorithm.

\begin{figure}[ht]
    \centering
    \begin{tikzpicture}	
        \draw (0,0) -- (3,0) -- (0,3) -- cycle;
        \draw[fill= gray!60!white] (1,0) -- (2,0) -- (2,1) -- (1,2) -- cycle;
        \draw [line width=0.7mm](0,0) -- (1,2);
        \draw [line width=0.7mm](1,0) -- (1,2);
        \draw [line width=0.7mm](2,0) -- (2,1);
        \draw [line width=0.7mm](0,3) -- (0,0);
        \draw [line width=0.7mm](1,0) -- (3,0);
        \draw [line width=0.7mm](1,2) -- (2,1);
        \draw node at (0,3) [regular_node]{};
        \draw node at (0,0) [regular_node]{};
        \draw node at (1,2) [regular_node]{};
        \draw node at (1,0) [regular_node]{};
        \draw node at (2,1) [regular_node]{};
        \draw node at (2,0) [regular_node]{};
        \draw node at (3,0) [regular_node]{};
        \node at (-0.7,1.5) {(1,1,1)};
        \node at (0.3,2) {(3)};
        \node at (0.6,0.5) {(2)};
        \node at (1.5,-0.5) {\small{type 1}};
    \end{tikzpicture}\hspace{0.5cm}
    \begin{tikzpicture}	
        \draw (0,0) -- (3,0) -- (0,3) -- cycle;
        \draw[fill= gray!60!white] (1,0) -- (2,0) -- (2,1) -- (1,1) -- cycle;
        \draw[fill= gray!60!white] (0,1) -- (1,2) -- (0,3) -- cycle;
        \draw [line width=0.7mm](0,1) -- (1,2);
        \draw [line width=0.7mm](0,3) -- (1,2);
        \draw [line width=0.7mm](0,3) -- (0,1);
        \draw [line width=0.7mm](1,0) -- (1,2);
        \draw [line width=0.7mm](2,1) -- (2,0);
        \draw [line width=0.7mm](1,0) -- (3,0);
        \draw [line width=0.7mm](1,1) -- (2,1);
        \draw (1,1) -- (0,0);
        \draw node at (0,3) [regular_node]{};
        \draw node at (0,1) [regular_node]{};
        \draw node at (1,2) [regular_node]{};
        \draw node at (1,1) [regular_node]{};
        \draw node at (1,0) [regular_node]{};
        \draw node at (2,1) [regular_node]{};
        \draw node at (2,0) [regular_node]{};
        \draw node at (3,0) [regular_node]{};
        \node at (-0.6,2) {(1,1)};
        \node at (1.5,-0.5) {\small{type 2}};
    \end{tikzpicture}
    \hspace{0.5cm}
    \begin{tikzpicture}	
        \draw (0,0) -- (3,0) -- (0,3) -- cycle;
        \draw[fill= gray!60!white] (1,0) -- (2,0) -- (2,1) -- (1,1) -- cycle;
        \draw[fill= gray!60!white] (0,1) -- (1,2) -- (0,3) -- cycle;
        \draw [line width=0.7mm](0,1) -- (1,2);
        \draw [line width=0.7mm](0,3) -- (1,2);
        \draw [line width=0.7mm](1,0) -- (1,2);
        \draw [line width=0.7mm](2,1) -- (2,0);
        \draw [line width=0.7mm](1,0) -- (3,0);
        \draw [line width=0.7mm](1,1) -- (2,1);
        \draw [line width=0.7mm](0,3) -- (0,1);
        \draw (1,1) -- (0,1);
        \draw node at (0,3) [regular_node]{};
        \draw node at (0,1) [regular_node]{};
        \draw node at (1,2) [regular_node]{};
        \draw node at (1,1) [regular_node]{};
        \draw node at (1,0) [regular_node]{};
        \draw node at (2,1) [regular_node]{};
        \draw node at (2,0) [regular_node]{};
        \draw node at (3,0) [regular_node]{};
        \node at (-0.6,2) {(1,1)};
        \node at (1.5,-0.5) {\small{type 3}};
    \end{tikzpicture}
    \begin{tikzpicture}	
        \draw (0,0) -- (3,0) -- (0,3) -- cycle;
        \draw[fill= gray!60!white] (0,0) -- (1,1) -- (2,0) -- cycle;
        \draw [line width=0.7mm](0,0) -- (0,3);
        \draw [line width=0.7mm](0,0) -- (3,0);
        \draw [line width=0.7mm](0,0) -- (1,1);
        \draw [line width=0.7mm](1,1) -- (2,0);
        \draw (1,1) -- (0,1);
        \draw (1,1) -- (1,2);
        \draw (2,1) -- (2,0);
        \draw node at (0,3) [regular_node]{};
        \draw node at (0,1) [regular_node]{};
        \draw node at (0,0) [regular_node]{};
        \draw node at (1,1) [regular_node]{};
        \draw node at (2,0) [regular_node]{};
        \draw node at (3,0) [regular_node]{};
        \node at (-0.6,2) {(1,1)};
        \node at (0.5,2) {(1,1)};
        \node at (1,-0.3) {(1,1)};
        \node at (1.5,-0.8) {\small{type 4}};
    \end{tikzpicture}
    \hspace{0.5cm}
    \begin{tikzpicture}	
        \draw (0,0) -- (3,0) -- (0,3) -- cycle;
        \draw [line width=0.7mm](0,0) -- (0,3);
        \draw [line width=0.7mm](0,0) -- (3,0);
        \draw (1,1) -- (2,0);
        \draw (1,1) -- (0,1);
        \draw (1,1) -- (1,2);
        \draw (2,1) -- (2,0);
        \draw node at (0,3) [regular_node]{};
        \draw node at (0,1) [regular_node]{};
        \draw node at (0,0) [regular_node]{};
        \draw node at (2,0) [regular_node]{};
        \draw node at (3,0) [regular_node]{};
        \node at (-0.6,2) {(1,1)};
        \node at (0.5,2) {(1,1)};
        \node at (1,0.2) {(1,1)};
        \node at (1,-0.3) {(1,1)};
        \node at (1.5,-0.8) {\small{type 5}};
    \end{tikzpicture}
    \hspace{0.5cm}
    \begin{tikzpicture}	
        \draw (0,0) -- (3,0) -- (0,3) -- cycle;
        \draw[fill= gray!60!white] (1,0) -- (2,0) -- (2,1) -- (1,1) -- cycle;
        \draw [line width=0.7mm](0,3) -- (0,1);
        \draw [line width=0.7mm](0,1) -- (2,1);
        \draw [line width=0.7mm](1,0) -- (3,0);
        \draw [line width=0.7mm](2,1) -- (2,0);
        \draw [line width=0.7mm](1,0) -- (3,0);
        \draw [line width=0.7mm](1,1) -- (1,0);
        \draw (1,1) -- (1,2);
        \draw node at (0,3) [regular_node]{};
        \draw node at (0,1) [regular_node]{};
        \draw node at (1,1) [regular_node]{};
        \draw node at (1,0) [regular_node]{};
        \draw node at (2,1) [regular_node]{};
        \draw node at (2,0) [regular_node]{};
        \draw node at (3,0) [regular_node]{};
        \node at (-0.6,2) {(1,1)};
        \node at (0.5,2) {(1,1)};
        \node at (1.5,-0.5) {\small{type 6}};
    \end{tikzpicture}
    \caption{Lattice path algorithm}\label{fig6}
\end{figure}
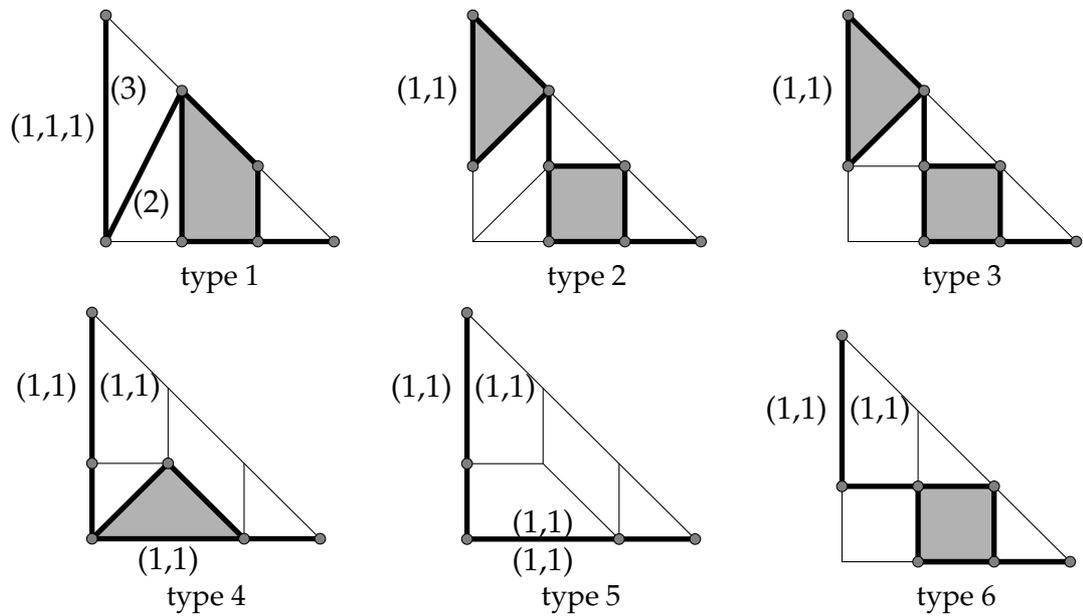

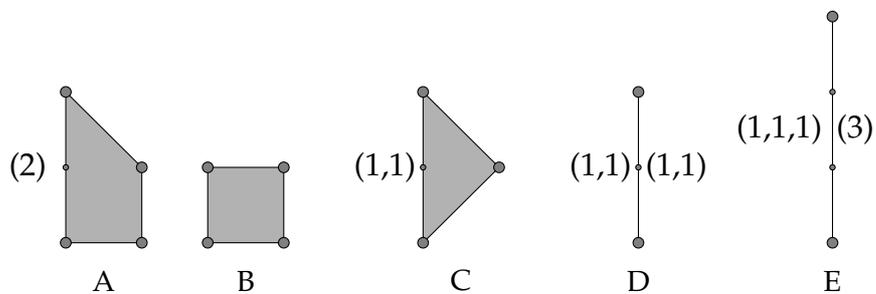
\begin{figure}[ht]	
	\centering
	\begin{tikzpicture}	
        \draw[fill= gray!60!white] (0,0) -- (1,0) -- (1,1) -- (0,2) -- cycle;
        \draw node at (0,0) [regular_node]{};
        \draw node at (0,2) [regular_node]{};
        \draw node at (1,1) [regular_node]{};
        \draw node at (1,0) [regular_node]{};
        \draw node at (0,1) [small_node]{};
        \node at (-0.5,1) {(2)};
        \node at (0.5,-0.5) {\small{A}};
    \end{tikzpicture}
    \hspace{0.5cm}
    \begin{tikzpicture}
        \draw[fill= gray!60!white] (0,0) rectangle (1,1);
        \draw node at (0,0) [regular_node]{};
        \draw node at (0,1) [regular_node]{};
        \draw node at (1,1) [regular_node]{};
        \draw node at (1,0) [regular_node]{};
        \node at (0.5,-0.5) {\small{B}};
    \end{tikzpicture}
    \hspace{0.5cm}
    \begin{tikzpicture}
        \draw[fill= gray!60!white] (0,0) -- (1,1) -- (0,2) -- cycle;
        \draw node at (0,0) [regular_node]{};
        \draw node at (1,1) [regular_node]{};
        \draw node at (0,2) [regular_node]{};
        \draw node at (0,1) [small_node]{};
        \node at (-0.5,1) {(1,1)};
        \node at (0.5,-0.5) {\small{C}};
    \end{tikzpicture}
    \hspace{0.5cm}
    \begin{tikzpicture}
        \draw (0,0) -- (0,2);
        \draw node at (0,0) [regular_node]{};
        \draw node at (0,2) [regular_node]{};
        \draw node at (0,1) [small_node]{};
        \node at (0,-0.5) {\small{D}};
        \node at (0.5,1) {(1,1)};
        \node at (-0.5,1) {(1,1)};
    \end{tikzpicture}
    \begin{tikzpicture}
        \draw (0,0) -- (0,3);
        \draw node at (0,0) [regular_node]{};
        \draw node at (0,1) [small_node]{};
        \draw node at (0,2) [small_node]{};
        \draw node at (0,3) [regular_node]{};
        \node at (0,-0.5) {\small{E}};
        \node at (0.3,1.5) {(3)};
        \node at (-0.7,1.5) {(1,1,1)};
    \end{tikzpicture}
    \caption{Rag rug elements}\label{fig7}
\end{figure}

\smallskip
{\bf Example.}
Let us compute the refined descendant invariant for the case
$$\Del=\left\{3\times(1,1),3\times(0,-1),3\times(-1,0)
\right\},\quad\overline{n}=(2,0,2,0,0,...)\ ,$$
considered in \cite[Example 9.28]{MR}. Here, $n=4$, and we choose the sequence $(m_1,m_2,m_3,m_4)=
(4,2,4,2)$, and the functional $\lambda(x,y)=x-\xi y$ with $0<\xi\ll1$. These data define
$$P(\Del)=\conv((0,0),(3,0),(0,3)),\quad q_{\min}=(0,3),\quad q_{\max}=(3,0)\ .$$
It is shown in \cite[Example 9.28]{MR} that there are $11$ admissible subdivisions, and they are presented in Figure \ref{fig6}.
The fat lines designate the $\lambda$-paths $\gamma^+_0$ and $\gamma^-_0$; the meaning of labels in parentheses we illustrate by an example:
in figure marked ``type 1", the segment $[(0,3),(0,0)]$ lies both, on $\gamma^+_0$ and $\gamma^-_0$, while in $\gamma^+_0$ it covers just one segment
and in $\gamma^-_0$ three segments of length $1$. Each of the figures marked ``type 1", ``type 2", ``type 3" represents one admissible
subdivision, each of the figures marked ``type 4" and ``type 6" represents two admissible subdivisions obtained by cutting the trapeze into
a triangle and a parallelogram, and, finally, figure marked ``type 5" represents four admissible subdivisions obtained by cutting each of the two trapezes into a triangle and a parallelogram (in fact, for each of these three types, different admissible subdivisions correspond to the
same isomorphism class of plane marked rational tropical curves).
Various type of rag rugs elements occurring in the admissible subdivisions are shown in Figures \ref{fig7}.

\medskip

The refined multiplicities of the admissible subdivision are (for types 4, 5, and 6 we sum up the multiplicities
over all obtained subdivisions):

\begin{equation*}
    \begin{aligned}
        \text{(1): }&\theta_y(\Del_E)\cdot [3]^{-}_y\cdot [2]^{-}_y\cdot [1]^{-}_y\cdot\theta_y(\Del_A)\cdot \frac{1}{3!}=
               16\cdot \frac{y^2+2y+3+2y^{-1}+y^{-2}}{y+2+y^{-1}}\\
        \text{(2): }&\theta_y(\Del_B)\cdot\theta_y(\Del_C) \cdot ([1]^{-}_y)^3 \cdot \frac {1}{2!}=
                16\cdot\frac{y^2+5y+6+5y^{-1}+y^{-2}}{(y+2+y^{-1})^2}\\
        \text{(3): }&\theta_y(\Del_B)\cdot\theta_y(\Del_C) \cdot ([1]^{-}_y)^3 \cdot \frac {1}{2!}=
                16\cdot\frac{y^2+5y+6+5y^{-1}+y^{-2}}{(y+2+y^{-1})^2}\\
        \text{(4): }&\theta_y(\Del_D)\cdot\theta_y(\Del_C) \cdot ([1]^{-}_y)^3 \cdot \frac {1}{(2!)^2}=
               48\cdot\frac{y+1+y^{-1}}{(y+2+y^{-1})^2}\\
        \text{(5): }&(\theta_y(\Del_D))^2 \cdot ([1]^{-}_y)^3 \cdot \frac {1}{(2!)^2}=
                \frac{144}{(y+2+y^{-1})^2}\\
        \text{(6): }&\theta_y(\Del_D)\cdot\theta_y(\Del_B) \cdot ([1]^{-}_y)^3 \cdot \frac {1}{2!}=
                48\cdot\frac{y+4+y^{-1}}{(y+2+y^{-1})^2}
    \end{aligned}
\end{equation*}
Finally, $$RD_y(\Del,\on)=
16\cdot\frac{y^3+6y^2+24y+46+24y^{-1}+6y^{-2}+y^{-3}}{(y+2+y^{-1})^2}$$
and
$$NRD_y(\Del,\on)=\frac{1}{36}RD_y(\Del,\on)=
\frac{4}{9}\cdot\frac{y^3+6y^2+24y+46+24y^{-1}+6y^{-2}+y^{-3}}{(y+2+y^{-1})^2}$$

Notice, that if we substitute $y=1$ into the latter expression, we get $NRD_1(\Del,\on)=3$
in agreement with \cite[Example 9.28]{MR}.

\smallskip

{\bf Remark.}
A possible generalization of the floor diagram algorithm as in \cite{BGM} seems to be more involved,
since the refined multiplicities of the floors do not admit reasonable explicit formulas contrary to the case
considered in \cite{BGM} and corresponding to $y=1$.

\medskip

{\bf Acknowledgements}. The authors have been supported by the German-Israeli Foundation
grant no. 1174-197.6/2011 and by the Israel Science Foundation grants no. 176/15 and 501/18,
as well as by the Bauer-Neuman Chair in Real and Complex Geometry.
This work has been started during the stay of the second author at the
Max-Planck Institut f\"ur Mathematik, Bonn, in August-September 2015, and then completed during the stay of the second author in
the Institute Mittag-Leffler, Stockholm, and \'Ecole Normale Sup\'erieure, Paris, in 2018. The second author is very grateful to
MPIM, IML, and ENS for hospitality and excellent working conditions. We also would like to thank Franziska Schroeter
for several important remarks and Travis Mandel for attracting our attention to the work \cite{M}.
Special thanks are due to the unknown referee for a careful reading of the paper and making many important critical remarks and suggestions.

{\ncsc School of Mathematical Sciences \\[-21pt]

Raymond and Beverly Sackler Faculty of Exact Sciences\\[-21pt]

Tel Aviv University \\[-21pt]

Ramat Aviv, 69978 Tel Aviv, Israel} \\[-21pt]

{\it E-mail address}: {\ntt blechm@gmail.com}, {\ntt shustin@tauex.tau.ac.il}

\end{document}